\documentclass[11pt]{article}
\usepackage{amsmath,amsthm,amssymb}
\usepackage[hyperfootnotes=false]{hyperref}
\usepackage{color}

\usepackage{titlesec}
\titleformat{\paragraph}
{\normalfont\normalsize\bfseries}{\theparagraph}{1em}{}
\titlespacing*{\paragraph}
{0pt}{3.25ex plus 1ex minus .2ex}{1.5ex plus .2ex}

\def\titlerunning#1{\gdef\titrun{#1}}
\makeatletter
\def\author#1{\gdef\autrun{\def\and{\unskip, }#1}\gdef\@author{#1}}
\def\address#1{{\def\and{\\\hspace*{18pt}}\renewcommand{\thefootnote}{}%
\footnote {#1}}%
\markboth{\autrun}{\titrun}}
\makeatother
\def\email#1{\hspace*{4pt}{\em e-mail}: #1}
\def\MSC#1{{\renewcommand{\thefootnote}{}%
\footnote{\emph{Mathematics Subject Classification (2010):} #1}}}
\def\keywords#1{\par\medskip
\noindent\textbf{Keywords:} #1}

\newtheorem{theorem}{Theorem}[section]

\newtheorem{cor}[theorem]{Corollary}
\newtheorem{lemma}[theorem]{Lemma}

\theoremstyle{definition}

\newtheorem{prob}[theorem]{Problem}
\newtheorem{cons}[theorem]{Construction}
\newtheorem{remark}[theorem]{Remark}

\numberwithin{equation}{section}

\def\bP{\mathbf P}

\def\bS{\mathbf S}

\def\cC{\mathcal C}
\def\cD{\mathcal D}

\def\cF{\mathcal F}
\def\cG{\mathcal G}
\def\cH{\mathcal H}

\def\cL{\mathcal L}
\def\cM{\mathcal M}

\def\cO{\mathcal O}
\def\cP{\mathcal P}
\def\cQ{\mathcal Q}

\def\cR{\mathcal R}
\def\cS{\mathcal S}

\def\cV{\mathcal V}
\def\cW{\mathcal W}
\def\cX{\mathcal X}
\def\cY{\mathcal Y}

\def\PG{{\rm PG}}

\def\GF{{\rm GF}}

\def\PGL{{\rm PGL}}

\def\PSp{{\rm PSp}}
\def\GL{{\rm GL}}
\def\SL{{\rm SL}}

\def\PGU{{\rm PGU}}

\def\rk{{\rm rk}}

\begin{document}


\baselineskip=16pt

\titlerunning{}

\title{On symmetric and Hermitian rank distance codes}

\author{Antonio Cossidente
\and Giuseppe Marino \and
Francesco Pavese}

\date{}

\maketitle

\address{A. Cossidente: Dipartimento di Matematica, Informatica ed Economia, Universit\`a degli Studi della Basilicata, Contrada Macchia Romana, 85100 Potenza, Italy; 
\email{antonio.cossidente@unibas.it}
\and
Giuseppe Marino: Dipartimento di Matematica e Applicazioni ``Renato Caccioppoli'', Universit\`a degli Studi di Napoli ``Federico II'', Complesso Universitario di Monte Sant'Angelo, Cupa Nuova Cintia, 21 - 80126 - Napoli, Italy; 
\email{giuseppe.marino@unina.it}
\and
F. Pavese: Dipartimento di Meccanica, Matematica e Management, Politecnico di Bari, Via Orabona 4, 70125 Bari, Italy; 
\email{francesco.pavese@poliba.it}}


\MSC{Primary 51E20; Secondary 05C70; 05A05}

\begin{abstract}
Let $\cM$ denote the set $\cS_{n, q}$ of $n \times n$ symmetric matrices with entries in $\GF(q)$ or the set $\cH_{n, q^2}$ of $n \times n$ Hermitian matrices whose elements are in $\GF(q^2)$.
Then $\cM$ equipped with the {\em rank distance} $d_r$ is a metric space. We investigate $d$--codes in $(\cM, d_r)$ and construct $d$--codes whose sizes are larger than the corresponding additive bounds. In the Hermitian case, we show the existence of an $n$--code of $\cM$, $n$ even and $n/2$ odd, of size $\left(3q^{n}-q^{n/2}\right)/2$, and of a $2$--code of size $q^6+ q(q-1)(q^4+q^2+1)/2$, for $n = 3$. In the symmetric case, if $n$ is odd or if $n$ and $q$ are both even, we provide better upper bound on the size of a $2$--code. In the case when $n = 3$ and $q>2$, a $2$--code of size $q^4+q^3+1$ is exhibited. This provides the first infinite family of $2$--codes of symmetric matrices whose size is larger than the largest possible additive $2$--code and an answer to a question posed in \cite[Section 7]{Sheekey}, see also \cite[p. 176]{S2}. 

\keywords{symmetric rank distance codes; Hermitian rank distance codes; symplectic polar spaces; Hermitian polar spaces; Segre variety.}
\end{abstract}

\section{Introduction}

Let $q$ be a power of a prime and let $\GF(q)$ be the finite field with $q$ elements. Denote by $\cS_{n, q}$ the set of $n \times n$ symmetric matrices with entries in $\GF(q)$ and by $\cH_{n, q^2}$ the set of $n \times n$ Hermitian matrices whose elements are in $\GF(q^2)$. Let $\cM$ be $\cS_{n, q}$ or $\cH_{n, q^2}$. For two matrices $A, B \in \cM$, define their {\em rank distance} to be 
$$
d_r(A, B) = \rk(A - B). 
$$
Thus $d_r$ is a metric on $\cM$ and $(\cM, d_r)$ is a metric space. A {\em rank metric code} $\cC$ is a non--empty subset of $(\cM, d_r)$. The {\em minimum distance} of $\cC$ is 
$$
d_r(\cC) = \min\{d_r(c_1,c_2) \;\; | \;\; c_1, c_2 \in \cC, c_1 \ne c_2\}. 
$$
We will refer to a code in $(\cM, d_r)$ with minimum distance $d$ as a {\em $d$--code}. A $d$--code is said to be {\em maximal} if it maximal with respect to set theoretic inclusion, whereas it is called {\em maximum} if it has the largest possible size. If a $d$--code $\cC \subset \cM$ forms a subgroup of $(\cM, +)$, then $\cC$ is called {\em additive}. Upper bounds on the size of a $d$--code of $\cM$ were provided in \cite[Corollary 7]{S0}, \cite[Lemma 3.5, Proposition 3.7]{S1}, \cite[Proposition 3.4]{S} and \cite[Theorems 1 and 2]{S3}. In the case when $\cC$ is additive much better bounds can be obtained. Indeed in \cite[Lemmas 3.5 and 3.6]{S1}, \cite[Theorem 4.3]{S2}, the author proved that the largest additive $d$--codes of $\cS_{n, q}$ have size at most either $q^{n(n-d+2)/2}$ or $q^{(n+1)(n-d+1)/2}$, according as $n-d$ is even or odd, respectively, whereas the size of the largest additive $d$--codes of $\cH_{n, q^2}$ cannot exceed $q^{n(n-d+1)}$, see \cite[Theorem 1]{S3}. Moreover there exist additive $d$--codes whose sizes meet the upper bounds for all possible value of $n$ and $d$, except when $\cM = \cH_{n, q^2}$, $n$, $d$ are both even and $3 < d < n$, see \cite[Theorems 12 and 16]{S0}, \cite[Theorem 4.4]{S1}, \cite[Theorem 5.3]{S2}, \cite[Theorems 4 and 5]{S3},\cite[Theorem 6.1]{S}, \cite{DGS}, \cite{GP}, \cite{GP1}. If $d$ is odd, a $d$--code attaining the corresponding additive bound is maximum. This is not always true if $d$ is even. However not much is known about $d$--codes whose size is larger than the corresponding additive bound. In the Hermitian case, if $n$ is even, there is an $n$--code of $\cH_{n, q^2}$ of size $q^n+1$ \cite[Theorem 6]{S3}, \cite[Theorem 18]{GLSV}. In the symmetric case only sporadic examples of non--additive $d$-codes that are larger than the largest possible additive $d$--code are known \cite[Tables 2 and 9]{S}.

\bigskip

Let $\cW(2n-1, q)$ be a non--degenerate symplectic polar space and $\cH(2n-1, q^2)$ be a non--degenerate Hermitian polar space. Let $\Pi_1$ be a generator of $\cW(2n-1, q)$ and let $\Lambda_1$ be a generator of $\cH(2n-1, q^2)$. It is known that there exists a bijection $\tau$ between the matrices of $\cS_{n, q}$ or $\cH_{n, q^2}$ and the generators of $\cW(2n-1, q)$ or $\cH(2n-1, q^2)$ disjoint from $\Pi_1$ or $\Lambda_1$, respectively, see \cite[Proposition 9.5.10]{BCN}, \cite{GLSV}. Here an upper bound on the maximum number of generators of $\cW(2n-1, q)$ or $\cH(2n-1, q^2)$ pairwise intersecting in at most an $(n-3)$--dimensional projective space is derived, see Theorems \ref{upper_symp} and \ref{upper_her}. As a by product, by means of $\tau$, the following upper bound on the size of a $2$--code $\cC$ of $\cS_{n, q}$ is obtained:
\begin{equation} \label{bound}
|\cC| \le 
\begin{cases}
\sum_{j = 0}^{\frac{n-1}{2}} \frac{q^j(q^{n-2j+1} + 1)}{q^{n-j+1} + 1} \prod_{i=1}^{j} \frac{q^{2(n-i+1)} - 1}{(q^i - 1)(q^{i-1} + 1)} & \mbox{ if } n \mbox{ is odd, } \\ 
\prod_{i = 2}^{n} (q^i + 1) & \mbox{ if } n \mbox{ is even. }
\end{cases}
\end{equation}
Since the previous known upper bound for the size of a $2$--code of $\cS_{n, q}$ was $q^{n(n-1)/2+1}(q^{n-1}+1)/(q+1)$ for $q$ odd \cite[Proposition 3.7]{S1}, and $q^{n(n+1)/2} - q^n + 1$ for $q$ even \cite[Proposition 3.4]{S}, it follows that \eqref{bound} provides better upper bounds if $n$ is odd or if $n$ and $q$ are both even. 

\bigskip

By using $\tau$, it can be seen that an $n$--code $\cC$ of $\cS_{n, q}$ or $\cH_{n, q^2}$ exists if and only if there exists a partial spread of $\cW(2n-1, q)$ or $\cH(2n-1, q^2)$ of size $|\cC|+1$, see Lemmas \ref{symplectic} and \ref{Hermitian}, \cite{GLSV}. It is well--known that the points of $\cW(2n -1, q)$ can be partitioned into $q^n+1$ pairwise disjoint generators of $\cW(2n-1, q)$, that is, $\cW(2n-1, q)$ admits a spread. On the other hand, $\cH(2n-1, q^2)$ has no spread. If $n$ is odd, an upper bound for the largest partial spreads of $\cH(2n-1, q^2)$ is $q^n+1$ \cite{V} and there are examples of partial spreads of that size \cite{L}. If $n$ is even the situation is less clear: upper bounds can be found in \cite{ISX}, as for lower bounds there is a partial spread of $\cH(2n-1, q^2)$ of size $(3q^2-q)/2 + 1$ for $n = 2$, $q > 13$, \cite[p. 32]{ACE} and of size $q^n+2$ for $n \ge 4$ \cite{GLSV}. Here, generalizing the partial spread of $\cH(3, q^2)$, we show the existence of a partial spread of $\cH(2n-1, q^2)$, in the case when $n$ is even and $n/2$ is odd, of size $\left(3q^{n}-q^{n/2}\right)/2 + 1$ (cf. Theorem \ref{partial_spread}) and hence, if $n$ is even and $n/2$ is odd, of an $n$--code of $\cH_{n, q^2}$ of size $\left(3q^{n}-q^{n/2}\right)/2$.

\bigskip

In the remaining part of the paper, we focus on the case $n = 3$. First, a further improvement on the size of a $2$--code $\cC$ of $\cS_{3, q}$ is obtained (cf. Corollary \ref{upper_sym3}):
$$
|\cC| \le \frac{q(q^2-1)(q^2+q+1)}{2}+1.  
$$
Then we construct $2$--codes of $\cS_{3, q}$ and $\cH_{3, q^2}$ of size $q^4+q^3+1$ and $q^6+ q(q-1)(q^4+q^2+1)/2$, respectively. This provides the first infinite family of $2$--codes of $\cS_{3, q}$ whose size is larger than the largest possible additive $2$--code and an answer to a question posed in \cite[Section 7]{Sheekey}, see also \cite[p. 176]{S2}. 

\section{Preliminaries}

\subsection{Projective and polar spaces}

Let $\PG(r - 1, q)$ be the projective space of projective dimension $r - 1$ over $\GF(q)$ equipped with homogeneous projective coordinates $X_1, \dots, X_r$. We will use the term $n$--space of $\PG(r - 1, q)$ to denote an $n$--dimensional projective subspace of $\PG(r - 1, q)$. We shall find it helpful to represent projectivities of $\PG(r - 1, q)$ by invertible $r \times r$ matrices over $\GF(q)$ and to consider the points of $\PG(r - 1, q)$ as column vectors, with matrices acting on the left. Let $U_i$ be the points having $1$ in the $i$--th position and $0$ elsewhere. Furthermore, we denote by $0_n$ and $I_n$ the $n \times n$ zero matrix and identity matrix, respectively; if $M$ is an $n \times n$ matrix over $\GF(q)$, we denote by $L(M)$ the $(n-1)$-space of $\PG(2n-1, q)$ whose underlying vector space is the vector space spanned by the rows of the $n \times 2n$ matrix $\begin{pmatrix} I_n & M \end{pmatrix}$; we also use the notation $L(M) = \left\langle \begin{pmatrix} I_n & M \end{pmatrix} \right\rangle$. If $m$ divides $r$, an {\em $(m-1)$--spread of $\PG(r-1, q)$} is a set of pairwise disjoint $(m-1)$--spaces of $\PG(r-1, q)$ which partition the point set of $\PG(r-1, q)$.

A {\em finite classical polar space} $\bP$ arises from a vector space of finite dimension over a finite field equipped with a non--degenerate reflexive sesquilinear form. In this paper we will be mainly concerned with symplectic polar spaces and Hermitian polar spaces. A projective subspace of maximal dimension contained in $\bP$ is called a {\em generator} of $\bP$. For further details on finite classical polar spaces we refer the readers to \cite{HT}. A {\em partial spread} $\bS$ of $\bP$ is a set of pairwise disjoint generators. A partial spread $\bS$ of $\bP$ is called a {\em spread of $\bP$} if $\bS$ partitions the point set of $\bP$.

\subsubsection{Segre varieties}
	
Consider the map defined by 	
$$
\xi: \PG(1,q)\times\PG(2, q) \longrightarrow \PG(5, q),
$$
taking a pair of points $x=(x_1, x_2)$ of $\PG(1,q)$, $y= (y_1, y_2, y_{3})$ of $\PG(2, q)$ to their product $(x_1y_1, \dots, x_2y_{3})$. This is a special case of a wider class of maps called {\em Segre maps} \cite{HT}.
The image of $\xi$ is an algebraic variety called the {\em Segre variety} and denoted by ${\Sigma}_{1,2}$. The Segre variety ${\Sigma}_{1,2}$ has two rulings, say ${\cal R}_1$ and ${\cal R}_2$, containing $q^2+q+1$ lines and $q+1$ planes, respectively, satisfying the following properties:  two subspaces in the same ruling are disjoint, elements of different ruling intersect in exactly one point and each point of ${\Sigma}_{1,2}$ is contained in exactly one member of each ruling. 

Notice that the set $\cR_1$ consists of all the lines of $\PG(5, q)$ incident with three distinct members of $\cR_2$ and, from \cite[Theorem 25.6.1]{HT}, three mutually disjoint planes of $\PG(5, q)$ define a unique Segre variety ${\Sigma}_{1, 2}$. A line of $\PG(5, q)$ shares $0,1,2$ or $q+1$ points with $\Sigma_{1,2}$. Also, the automorphism group of ${\Sigma}_{1,2}$ in $\PGL(6, q)$ is a group isomorphic to $\PGL(2, q) \times \PGL(3, q)$ \cite[Theorem 25.5.13]{HT}. For more details on Segre varieties, see \cite{HT}.

\subsection{Graphs}

Recall some definitions and results from \cite{BH, GR}. Suppose $\Gamma$ is a (simple, undirected) graph having $\cV$ as set of vertices. The adjacency matrix $A$ of $\Gamma$ is a symmetric real matrix whose rows and columns are indexed by $1, \dots, |\cV|$. The eigenvalues of $\Gamma$ are those of its adjacency matrix $A$. A graph $\Gamma$ is called {\em regular of valency} $k$ or {\em $k$--regular} when every vertex has precisely $k$ neighbors. If $\Gamma$ is regular of valency $k$, then $A {\bf 1} = k {\bf 1}$, where ${\bf 1}$ denotes the all one column vector. Hence $k$ is an eigenvalue of $\Gamma$ and for every eigenvalue $\lambda$ of $\Gamma$, we have that $|\lambda| \le k$. Furthermore the multiplicity of $k$ equals the number of connected components of $\Gamma$. 

Let $\Gamma$ be a $k$--regular graph and let $\{\cV_1, \dots, \cV_m\}$ be a partition of $\cV$. Let $A$ be partitioned according to $\{\cV_1, \dots, \cV_m\}$, that is,
		\[A = \begin{pmatrix}
		A_{1,1} & \dots & A_{1,m} \\
		\vdots & & \vdots \\
		A_{m,1} & \dots & A_{m,m}
		\end{pmatrix},\]
such that $A_{i, i}$ is a square matrix for all $1 \leq i \leq m$. The {\em quotient matrix} $B$ is the $m \times m$ matrix with entries the average row sum of the blocks of $A$. More precisely, 
		$$
		B = (b_{i, j}),  \; b_{i, j} = \frac{1}{v_i} {\bf 1}^t A_{i, j} {\bf 1} ,
		$$
where $v_i$ is the number of rows of $A_{i, j}$. If the row sum of each block $A_{i, j}$ is constant then the partition is called {\em equitable} or
{\em regular} and we have $A_{i, j}{\bf 1} = b_{i, j}{\bf 1}$ for $1 \le i, j \le m$. One important class of equitable partitions arises from automorphisms of $\Gamma$, indeed the orbits of any group of automorphisms of $\Gamma$ form an equitable partition. The following result is well known and useful. 

\begin{lemma}[Lemma 2.3.1, \cite{BH}, Theorem 9.4.1, \cite{GR}]\label{equitable}
Let $B$ be the quotient matrix of an equitable partition. If $\lambda$ is an eigenvalue of $B$, then $\lambda$ is an eigenvalue of $A$.

Let $\Gamma$ be a vertex--transitive graph and let $B$ be the quotient matrix of an equitable partition arising from the orbits of some subgroup of $Aut(\Gamma)$. If $|\cV_i| = 1$ for some $i$, then every eigenvalue of $A$ is an eigenvalue of $B$.
\end{lemma}		

A {\em coclique} of $\Gamma$ is a set of pairwise nonadjacent vertices. The independence number $\alpha(\Gamma)$ is the size of the largest coclique of $\Gamma$. Let $\lambda_1 \ge \dots \ge \lambda_{|\cV|}$ be the eigenvalues of $\Gamma$. The following results are due to Cvetkovi{\'c} \cite{C} and Hoffman, respectively; see also \cite[Theorem 3.5.1, Theorem 3.5.2]{BH}.

\begin{lemma}\label{cvetkovic}
$\alpha(\Gamma) \le \min\{|\{i \;\; | \;\; \lambda_i \ge 0\}|, |\{i \;\; | \;\; \lambda_i \le 0\}|\}$.
\end{lemma}

\begin{lemma}\label{hoffman}
$\alpha(\Gamma) \le -\frac{|\cV| \lambda_{|\cV|}}{k- \lambda_{|\cV|}}$.
\end{lemma}

\section{Symmetric matrices and symplectic polar spaces}

Let $\cW(2n-1, q)$ be the non--degenerate symplectic polar space of $\PG(2n-1, q)$ associated with the following alternating bilinear form 
$$
\begin{pmatrix} X_1, \dots, X_{2n} \end{pmatrix}
\begin{pmatrix}
0_n & I_n \\
- I_n & 0_n
\end{pmatrix}
\begin{pmatrix} Y_1 \\ \vdots \\ Y_{2n} \end{pmatrix}.
$$ 
Let $\perp$ denote the symplectic polarity of $\PG(2n-1, q)$ defining $\cW(2n-1, q)$ and let $\PSp(2n, q) = {\rm Sp}(2n, q)/\langle -I \rangle$, where ${\rm Sp}(2n, q)$ is the group of isometries of the alternating bilinear form previously defined. Hence $\PSp(2n, q)$ consists of projectivities of $\PG(2n-1, q)$ fixing $\cW(2n-1, q)$. It acts transitively on the generators of $\cW(2n-1, q)$. Denote by $\Pi_1$ the $(n-1)$--space of $\PG(2n -1, q)$ spanned by $U_{n+1}, \dots, U_{2n}$. Then $\Pi_1$ is a generator of $\cW(2n-1, q)$. Let $G$ be the stabilizer of $\Pi_1$ in $\PSp(2n, q)$. Then it is readily seen that an element of $G$ is represented by the matrix
\begin{equation}
\begin{pmatrix} \label{matrice_s}
T^{-t} & 0_n \\
S_0 T^{-t} & T \\
\end{pmatrix},
\end{equation}
where $T \in \GL(n, q)$ and $S_0 \in \cS_{n, q}$. Hence $G \simeq \cS_{n, q} \rtimes \left(\GL(n, q)/\langle -I_{2n} \rangle\right)$ has order 
$$
\frac{q^{n(n+1)/2} \prod_{i = 0}^{n-1} (q^n-q^i)}{gcd(2, q-1)}
$$
and it acts transitively on the set $\cG$ of generators of $\cW(2n-1, q)$ disjoint from $\Pi_1$.

Define an action of $\cS_{n, q} \rtimes \GL(n, q)$ on $\cS_{n, q}$ as follows
$$
\left((S_0, T), S \right) \in \left(\cS_{n, q} \rtimes \GL(n, q)\right) \times \cS_{n, q} \longmapsto T S T^t + S_0 \in \cS_{n, q}.
$$
Its orbitals are the relations of an association scheme, the so called association scheme of symmetric matrices \cite{HW, WJ}. The following result enlightens a correspondence between $\cS_{n, q}$ and $\cG$, see also \cite[Proposition 9.5.10]{BCN}.

\begin{lemma}\label{symplectic}
There is a bijection between $\cS_{n, q}$ and $\cG$ such that $\cS_{n, q} \rtimes \GL(n, q)$ acts on $\cS_{n, q}$ as $G$ acts on $\cG$. In particular, a $d$--code of $\cS_{n, q}$ corresponds to a set of generators of $\cW(2n-1, q)$ disjoint from $\Pi_1$ pairwise intersecting in at most an $(n - d - 1)$--space, and conversely.
\end{lemma}
\begin{proof}
Let $S \in \cS_{n, q}$. Since the rank of the matrix $\begin{pmatrix} I_n & S \\ 0_n & I_n \end{pmatrix}$ is $2n$, it follows that $L(S)$ is disjoint from $\Pi_1$. The map $S \mapsto L(S)$ is injective. Moreover $|\cS_{n, q}| = |\cG|$ and $L(S)$ is a generator of $\cW(2n-1, q)$, indeed 
$$
\begin{pmatrix} I_n & S \end{pmatrix}  
\begin{pmatrix}
0_n & I_n \\
- I_n & 0_n
\end{pmatrix}
\begin{pmatrix} I_n \\ S \end{pmatrix} = 0 .
$$
Finally, let $g \in G$ be represented by the matrix \eqref{matrice_s}. Then $L(S)^g = \left\langle \begin{pmatrix} T^{-1} & S T^t + T^{-1} S_0 \end{pmatrix} \right\rangle = \left\langle \begin{pmatrix} I_n & T S T^t + S_0 \end{pmatrix} \right\rangle = L(TST^t + S_0)$. This completes the proof of the first part of the statement. Let $\cC$ be a $d$-code of $\cS_{n,q}$ and let $S_1$ and $S_2$ be two different elements of $\cC$. Since $\rk (S_1-S_2) \geq d$, it follows that $L(S_1) \cap L(S_2)$ is at most an $(n-d-1)$--space.
\end{proof}
Let $\Pi_2 = L(0_n)$. The previous lemma implies that the number of orbits of $G_{\Pi_2}$ on $\cG$ equals the number of relations of the association scheme on symmetric matrices. Since $|\cG| = q^{n(n+1)/2}$ and $G$ acts transitively on $\cG$, it follows that the stabilizer of $\Pi_2$ in $G$, namely $G_{\Pi_2}$, has order
$$
\frac{\prod_{i = 0}^{n-1} (q^n-q^i)}{gcd(2, q-1)}.
$$
More precisely an element of $G_{\Pi_2}$ is represented by the matrix
\begin{equation}
\begin{pmatrix} \label{matrice_s1}
T^{-t} & 0_n \\
0_n & T \\
\end{pmatrix},
\end{equation}
where $T \in \GL(n, q)$. The group $G_{\Pi_2}$ acts transitively on points and hyperplanes of both $\Pi_1$ and $\Pi_2$. The action of the group $G_{\Pi_2}$ on points of $\cW(2n-1, q)$ has been studied in \cite[p. 347]{K}. For the sake of completeness a direct proof is given below. 

\begin{lemma}\label{point_orbits}
The orbits of $G_{\Pi_2}$ on points of $\cW(2n-1, q) \setminus (\Pi_1 \cup \Pi_2)$ are 
\begin{itemize}
\item $\cP_0$ of size $(q^n-1)(q^{n-1}-1)/(q-1)$,
\item $\cP_1$ of size $q^{n-1}(q^n-1)$,
\end{itemize}
if $q$ is even and
\begin{itemize}
\item $\cP_0$ of size $(q^n-1)(q^{n-1}-1)/(q-1)$,
\item $\cP_1, \cP_2$, both of size $q^{n-1}(q^n-1)/2$,
\end{itemize} 
if $q$ is odd.
\end{lemma}
\begin{proof}
Let $g$ be the projectivity of $G_{\Pi_2}$ associated with the matrix \eqref{matrice_s1}, for some $T \in \GL(n , q)$. Then $g$ stabilizes $U_{n+1}$ if and only if the first column of $T$ is $(z, 0, \dots, 0)^t$, for some $z \in \GF(q) \setminus \{0\}$, which is equivalent to the requirement that the first row of $T^{-t}$ is $(z^{-1}, 0, \dots, 0)$, for some $z \in \GF(q) \setminus \{0\}$. It follows that $Stab_{G_{\Pi_2}}(U_{n+1})$ has two orbits on points of $\Pi_2$, namely $U_{n+1}^\perp \cap \Pi_2$ and $\Pi_2 \setminus U_{n+1}^\perp$. Hence $G_{\Pi_2}$ permutes in a single orbit the lines of $\cW(2n-1, q)$ meeting both $\Pi_1$, $\Pi_2$ in a point. Similarly the lines not of $\cW(2n-1, q)$ meeting both $\Pi_1$, $\Pi_2$ in a point form a unique $G_{\Pi_2}$--orbit.

Let $P = U_2 + U_{n+1}$. A projectivity of $G_{\Pi_2}$ stabilizing $P$ has to fix both $U_2$ and $U_{n+1}$. Straightforward calculations show that a member of $G_{\Pi_2}$ fixes $P$ if and only if it is associated with the matrix \eqref{matrice_s1}, where
$$
T = \begin{pmatrix} 
x & * & * & \dots & * \\
0 & x^{-1} & 0 & \dots & 0 \\
0 & * & & & & \\
\vdots & \vdots & & T' & & \\
0 & * & & & & \\
\end{pmatrix},
$$
for some $T' \in \GL(n - 2, q)$ and $x \in \GF(q) \setminus \{0\}$. Therefore $|Stab_{G_{\Pi_2}}(P)| = (q-1)(q^{n-1}-1)(q^{n-2}-1) |\GL(n-2, q)|/gcd(2, q-1)$. Hence $|P^{G_{\Pi_2}}| = (q^n-1)(q^{n-1}-1)/(q-1)$, which equals the number of points of $\cW(2n-1, q) \setminus (\Pi_1 \cup \Pi_2)$ that lie on the lines of $\cW(2n-1, q)$ meeting both $\Pi_1$, $\Pi_2$ in one point.

Let $z \in \GF(q) \setminus \{0\}$ and let $P = U_1 + z U_{n+1}$. As before, a projectivity of $G_{\Pi_2}$ stabilizing $P_z$ has to fix both $U_1$ and $U_{n+1}$. Straightforward calculations show that the projectivity $g \in G_{\Pi_2}$ fixes the line $U_1 U_{n+1}$ if and only if it is associated with the matrix \eqref{matrice_s1}, where
$$
T = \begin{pmatrix} 
y & 0 & \dots & 0 \\
0 & & & \\
\vdots & & T'' & \\
0 &  &  & \\
\end{pmatrix},
$$
for some $T'' \in \GL(n - 1, q)$ and $y \in \GF(q) \setminus \{0\}$. Moreover $g$ stabilizes $P_z$ if and only if $y = \pm 1$. In this case we have that $|Stab_{G_{\Pi_2}}(P)| = |\GL(n-1, q)|$. Hence if $q$ is even, $|P_z^{G_{\Pi_2}}| = q^{n-1}(q^{n}-1)$, which equals the number of points of $\cW(2n-1, q) \setminus (\Pi_1 \cup \Pi_2)$ that lie on the lines not belonging to $\cW(2n-1, q)$ and meeting both $\Pi_1$, $\Pi_2$ in one point. If $q$ is odd, then $|P_z^{G_{\Pi_2}}| = q^{n-1}(q^{n}-1)/2$. Representatives for these two orbits are $P_{z_1}$ and $P_{z_2}$, where $z_1$ is a non--zero square in $\GF(q)$ and $z_2$ is a non--square of $\GF(q)$. Indeed there is no element of $G_{\Pi_2}$ sending $P_{z_1}$ to $P_{z_2}$. To see this fact assume on the contrary that there is a projectivity of $G_{\Pi_2}$ mapping $P_{z_1}$ to $P_{z_2}$. Then it has to fix the line $U_1 U_{n+1}$. On the other hand such a projectivity sends the point $P_{z_1}$ to $P_{y^2 z_1}$. Hence $z_2 = y^2 z_1$, a contradiction.   
\end{proof}

\begin{remark}\label{oss}
Note that if $q$ is odd and $\ell$ is a line such that $\ell \cap \Pi_2 = P_2 = (x_1, \dots, x_{n}, 0, \dots, 0)$, $\ell \cap \Pi_1 = P_1 = (0, \dots, 0, y_{n+1}, \dots, y_{2n})$ and $\ell$ is not a line of $\cW(2n-1, q)$, then the point $P = P_1 + z P_2 \in \ell$ belongs to $\cP_1$ or to $\cP_2$, according as $z (x_1, \dots, x_{n}) (y_{n+1}, \dots, y_{2n})^t$ is a non--zero square or a non--square in $\GF(q)$. Therefore $|\ell \cap \cP_1| = |\ell \cap \cP_2| = (q-1)/2$. Moreover, it can be checked that there are projectivities of $Stab_{\PSp(2n-1, q)}(\{\Pi_1, \Pi_2\}) \setminus G_{\Pi_2}$ interchanging the two orbits and hence $Stab_{\PSp(2n-1, q)}(\{\Pi_1, \Pi_2\})$ acts transitively on points of $\cW(2n-1, q) \setminus (\Pi_1 \cup \Pi_2)$.
\end{remark}

From Lemma \ref{symplectic}, $\Pi_3$ is a generator of $\cW(2n-1, q)$ disjoint from both $\Pi_1$ and $\Pi_2$ if and only if $\Pi_3 = L(A)$, where $A$ is an invertible matrix of $\cS_{n, q}$. Hence there are 
$$
q^{\frac{n(n+1)}{2} - \left\lceil \frac{n}{2} \right\rceil^2} \prod_{i = 1}^{\left\lceil \frac{n}{2} \right\rceil} \left( q^{2i-1}-1 \right)
$$
generators of $\cW(2n-1, q)$ disjoint from both $\Pi_1$ and $\Pi_2$, see \cite[Theorem 2]{Mac}, \cite[Corollary 19]{KMS}. The following result is well known. For the convenience of the reader a direct proof is provided.

\begin{lemma}[Theorem 21, \cite{KMS}]\label{segre_symp}
Let $\Pi_3$ be a generator of $\cW(2n-1, q)$ disjoint from $\Pi_1$ and $\Pi_2$. The points $P \in \Pi_3$ such that there exists a line of $\cW(2n-1, q)$ through $P$ intersecting $\Pi_1$ and $\Pi_2$ are the absolute points of a non--degenerate polarity which is
\begin{itemize}
\item pseudo--symplectic if $q$ is even and $n$ is odd,
\item orthogonal if $q$ and $n$ are odd,
\item symplectic or pseudo--symplectic if $q$ and $n$ are even,
\item elliptic orthogonal or hyperbolic orthogonal if $q$ is odd and $n$ is even. 
\end{itemize}
\end{lemma}
\begin{proof}
Let $\Pi_3 = L(A)$ be a generator of $\cW(2n-1, q)$ such that $|\Pi_1 \cap \Pi_3| = |\Pi_2 \cap \Pi_3| = 0$. We show that there is a non--degenerate polarity $\rho$ of $\Pi_3$ associated with the matrix $A$. Observe that the $(n-1)$--space $\Pi_3$ has equations:
$$
\begin{pmatrix}
X_{n+1} \\
\vdots \\
X_{2n}
\end{pmatrix} = A 
\begin{pmatrix}
X_1 \\
\vdots \\
X_n
\end{pmatrix}
$$
Hence the point $P$ belongs to $\Pi_3$ if and only if $P = (x_1, \dots, x_n, 0, \dots, 0) + (0, \dots, 0, x_1, \dots, x_n) \begin{pmatrix} 0_n & A \end{pmatrix}^t$ and $P$ lies on the line $\ell$ joining the points $(0, \dots, 0, x_1, \dots, x_n) \begin{pmatrix} 0_n & A \end{pmatrix}^t \in \Pi_1$ and $(x_1, \dots, x_n, 0, \dots, 0) \in \Pi_2$. Thus $\ell^\perp$ is represented by the equations: $x_1 X_{n+1} + \ldots + x_n X_{2n} = (x_1, \dots, x_n) A \left(X_1, \dots, X_n \right)^t = 0$ and a point $P' = (y_1, \dots, y_n, 0, \dots, 0) + (0, \dots, 0, y_1, \dots, y_n) \begin{pmatrix} 0_n & A \end{pmatrix}^t \in \Pi_3$ belongs to $\ell^\perp$ if and only if 
$$
(y_1, \dots, y_n) A 
\begin{pmatrix}
x_1 \\ 
\vdots \\ x_n
\end{pmatrix} = 0.
$$ 
This concludes the proof.
\end{proof}

\begin{lemma}\label{dis}
The group $G_{\Pi_2}$ has the following orbits on generators of $\cW(2n-1, q)$ disjoint from both $\Pi_1$ and $\Pi_2$:
\begin{itemize}
\item one orbit if $q$ is even and $n$ is odd,
\item two equally sized orbits if $q$ and $n$ are odd, 
\item two orbits having size 
$$
q^{\frac{n(n-2)}{4}} \prod_{i = 1}^{\frac{n}{2}} \left( q^{2i-1} - 1 \right) \quad \mbox{ and } \quad q^{\frac{n(n-2)}{4}} \left( q^n -1 \right) \prod_{i = 1}^{\frac{n}{2}} \left( q^{2i-1} - 1 \right)
$$ 
if $q$ and $n$ are even,
\item two orbits having size 
$$
\frac{q^{\frac{n^2}{4}} \left( q^{\frac{n}{2}} + 1 \right) \prod_{i = 1}^{\frac{n}{2}} \left( q^{2i-1} - 1 \right)}{2} \quad \mbox{ and } \quad \frac{q^{\frac{n^2}{4}} \left( q^{\frac{n}{2}} - 1 \right) \prod_{i = 1}^{\frac{n}{2}} \left( q^{2i-1} - 1 \right)}{2}
$$ 
if $q$ is odd and $n$ is even.   
\end{itemize}
\end{lemma}
\begin{proof}
Let $\Pi_3 = L(A)$ be a generator of $\cW(2n-1, q)$ such that $|\Pi_1 \cap \Pi_3| = |\Pi_2 \cap \Pi_3| = 0$ and let $g$ be the projectivity of $G_{\Pi_2}$ associated with the matrix \eqref{matrice_s1}, for some $T \in \GL(n , q)$. From the proof of Lemma \ref{point_orbits}, $\Pi_3^g = L(TAT^t)$. Therefore $g$ stabilizes $\Pi_3$ if and only if $TAT^t = A$. It follows that \cite[Appendix I]{H}
$$
|Stab_{G_{\Pi_2}}(\Pi_3)| = 
\begin{cases}
q^{\frac{(n-1)^2}{4}} \prod_{i = 1}^{\frac{n-1}{2}} \left( q^{2i} - 1 \right) & \mbox{ for } n \mbox{ odd,} \\

q^{\frac{n^2}{4}} \prod_{i = 1}^{\frac{n}{2}} \left( q^{2i} - 1 \right) & \mbox{ for } q, n \mbox{ even, } a_{i j} = a_{j i}, a_{i i} = 0, \\

q^{\frac{n^2}{4}} \prod_{i = 1}^{\frac{n - 2}{2}} \left( q^{2i} - 1 \right) & \mbox{ for } q, n \mbox{ even, } a_{i j} = a_{j i}, a_{i i} \ne 0, \mbox{ for some } i, \\

q^{\frac{n(n-2)}{4}} \left( q^{\frac{n}{2}} - 1 \right) \prod_{i = 1}^{\frac{n-2}{2}} \left( q^{2i} - 1 \right) & \mbox{ for } q \mbox{ odd, } n \mbox{ even, } \det(A) \mbox{ square of } \GF(q) \setminus \{0\}, \\

q^{\frac{n(n-2)}{4}} \left( q^{\frac{n}{2}} + 1 \right) \prod_{i = 1}^{\frac{n-2}{2}} \left( q^{2i} - 1 \right)& \mbox{ for } q \mbox{ odd, } n \mbox{ even, } \det(A) \mbox{ non--square of } \GF(q), \\
\end{cases}
$$
where $A = (a_{i j})$. The result follows.
\end{proof}

\begin{remark}\label{oss1}
We remark that if $q$ and $n$ are odd, then the generator $\Pi_3 = L(A)$, $A \in \cS_{n, q}$, $\rk(A) = n$, belongs to the first or the second $G_{\Pi_2}$--orbit on generators of $\cW(2n-1, q)$ skew to $\Pi_1$, $\Pi_2$, according as $\det(A)$ is a square or a non--square in $\GF(q)$. It can be easily seen that there are projectivities of $Stab_{\PSp(2n-1, q)}(\{\Pi_1, \Pi_2\}) \setminus G_{\Pi_2}$ interchanging the two orbits. Hence $Stab_{\PSp(2n-1, q)}(\{\Pi_1, \Pi_2\})$ acts transitively on generators of $\cW(2n-1, q)$ disjoint from $\Pi_1$ and $\Pi_2$.
\end{remark}

\subsection{$\cW(5, q)$}

Set $n = 3$. Let $\cW(5, q)$ be the symplectic polar space of $\PG(5, q)$ as described above. Recall that $\cG$ is the set of $q^6$ planes of $\cW(5, q)$ that are disjoint from $\Pi_1$, the group $G$ is the stabilizer of $\Pi_1$ in $\PSp(6, q)$, $\Pi_2 = L(0_3)$ and $G_{\Pi_2}$ is the stabilizer of $\Pi_2$ in $G$. 

Following Lemma \ref{point_orbits}, let $\cP_0$ be the set of points $R$ of $\cW(5, q) \setminus (\Pi_1 \cup \Pi_2)$ such that the line through $R$ intersecting $\Pi_1$ and $\Pi_2$ is a line of $\cW(5, q)$ and let $\cP$ be its complement in $\cW(5, q) \setminus (\Pi_1 \cup \Pi_2)$. Note that $\cP$ coincides with $\cP_1$ or $\cP_1 \cup \cP_2$, according as $q$ is even or odd. Then $|\cP_0| = (q^2-1) (q^2+q+1)$ and $|\cP| = q^5 - q^2$. Let $\ell$ be a line of $\cW(5, q)$ disjoint from $\Pi_1 \cup \Pi_2$. The hyperplane $\langle \Pi_2, \ell \rangle$ meets $\Pi_1$ in a line, say $r_\ell$, and the three--space $\langle \ell, r_{\ell} \rangle$ meets $\Pi_2$ in a line, say $t_{\ell}$. Hence the line $\ell$ defines a unique three--space $T_\ell = \langle r_\ell, t_\ell \rangle$ meeting both $\Pi_1$, $\Pi_2$ in a line. 

\begin{lemma}
Let $\ell$ be a line of $\cW(5, q)$ disjoint from $\Pi_1 \cup \Pi_2$, then $|\ell \cap \cP_0|$ belongs to $\{1, q+1\}$, if $q$ is even, and to $\{0, 1, 2\}$, if $q$ is odd.
\end{lemma}
\begin{proof}
There are two possibilities, either $T_\ell^\perp$ is a line of $\cW(5, q)$ or it is not. In the former case, among the $q+1$ lines meeting $\ell$, $r_\ell$, $t_\ell$ in one point, there is exactly one line of $\cW(5, q)$. If the latter case occurs, then $T_\ell \cap \cW(5, q)$ is a $\cW(3, q)$ and the regulus $\cR$ determined by $\ell, r_\ell, t_\ell$ consists of lines of $\cW(3, q)$. Thus its opposite regulus consists of either $1$ or $q+1$ lines of $\cW(3, q)$ if $q$ is even and of $0$ or $2$ lines of $\cW(3, q)$ if $q$ is odd.
\end{proof}

Let us partition the set of lines of $\cW(5, q)$ disjoint from $\Pi_1 \cup \Pi_2$. Let  
\begin{eqnarray*}
&& \cL_0 = \{\ell \mbox{ line of } \cW(5, q) \;\; : \;\; |\ell \cap (\Pi_1 \cup \Pi_2)| = 0, |\ell \cap \cP_0| = q+1 \},\\
&& \cL_1 = \{\ell \mbox{ line of } \cW(5, q) \;\; : \;\; |\ell \cap (\Pi_1 \cup \Pi_2)| = 0, T_\ell^\perp \mbox{ is a line of } \cW(5, q) \},\\
&& \cL_2 = \{\ell \mbox{ line of } \cW(5, q) \;\; : \;\; |\ell \cap (\Pi_1 \cup \Pi_2)| = 0, |\ell \cap \cP_0| = 1, T_\ell^\perp \mbox{ is not a line of } \cW(5, q) \},
\end{eqnarray*}
if $q$ is even, or
\begin{eqnarray*}
&& \cL_0 = \{\ell \mbox{ line of } \cW(5, q) \;\; : \;\; |\ell \cap (\Pi_1 \cup \Pi_2)| = 0, |\ell \cap \cP_0| =  0, T_\ell^\perp \mbox{ is not a line of } \cW(5, q) \},\\
&& \cL_1 = \{\ell \mbox{ line of } \cW(5, q) \;\; : \;\; |\ell \cap (\Pi_1 \cup \Pi_2)| = 0, T_\ell^\perp \mbox{ is a line of } \cW(5, q) \},\\
&& \cL_2 = \{\ell \mbox{ line of } \cW(5, q) \;\; : \;\; |\ell \cap (\Pi_1 \cup \Pi_2)| = 0, |\ell \cap \cP_0| =  2, T_\ell^\perp \mbox{ is not a line of } \cW(5, q) \},
\end{eqnarray*}
if $q$ is odd.
Note that in both cases if $\ell \in \cL_1$, then $|\ell \cap \cP_0| =1$, whereas if $q$ is even and $\ell \in \cL_0$, then $T_\ell^\perp$ is not a line of $\cW(5,q)$.

\begin{lemma}
If $q$ is even, then
$$
|\cL_0| = q^2(q^3-1), \quad |\cL_1| = q(q^2-1)(q^3-1), \quad |\cL_2| = q^2(q^2-1)(q^3-1).
$$
If $q$ is odd, then 
$$
|\cL_0| = \frac{q^3(q-1)(q^3-1)}{2}, \quad |\cL_1| = q(q^2-1)(q^3-1), \quad |\cL_2| = \frac{q^3(q+1)(q^3-1)}{2}.
$$
\end{lemma}
\begin{proof}
The line $T_\ell^\perp$ meets both $\Pi_1$, $\Pi_2$ in one point. If $T_\ell^\perp$ is a line of $\cW(5, q)$, then $T_\ell \cap \cW(5, q)$ consists of $q+1$ generators of $\cW(5, q)$ through $T_\ell^\perp$ and hence there are $q(q-1)^2$ lines of $\cW(5, q)$ contained in $T_\ell$ and disjoint from $\Pi_1 \cup \Pi_2$. Since there are $(q+1)(q^2+q+1)$ lines of $\cW(5, q)$ meeting both $\Pi_1$, $\Pi_2$ in one point, we get $|\cL_1| = q(q^2-1)(q^3-1)$. If $T_\ell^\perp$ is not a line of $\cW(5, q)$, then $T_\ell \cap \cW(5, q)$ is a non--degenerate symplectic polar space $\cW(3, q)$ and the regulus $\cR$ determined by $\ell, r_\ell, t_\ell$ is a regulus of $\cW(3, q)$. The point line dual of $\cW(3, q)$ is a parabolic quadric $\cQ(4, q)$ and the lines $r_\ell$ and $t_\ell$ correspond to two points $R, T$ such that the line $R T$ meets $\cQ(4, q)$ only in $R$ and $T$. Moreover, the regulus $\cR$ corresponds to a conic $C$ of $\cQ(4, q)$, where $R, T \in C$.

Assume that $q$ is even, then $\ell$ belongs either to $\cL_0$ or to $\cL_2$, according as the opposite regulus of $\cR$ has $q+1$ or one line of $\cW(3, q)$. In this case the parabolic quadric $\cQ(4, q)$ has a nucleus, say $N$. Moreover, the opposite regulus of $\cR$ has $q+1$ or one line of $\cW(3, q)$ according as $N$ belongs to the plane $\langle C \rangle$ or does not. Therefore, in $\cW(3, q)$, $\ell$ can be chosen in $q-1$ ways such that it belongs to $\cL_0$ and in $(q^2-1)(q-1)$ ways such that it belongs to $\cL_2$. Since there are $q^2(q^2+q+1)$ lines not of $\cW(5, q)$ meeting both $\Pi_1$, $\Pi_2$ in one point, we get $|\cL_0| = q^2(q^3-1)$ and $|\cL_2| = q^2(q^2-1)(q^3-1)$.   

If $q$ is odd, then $\ell$ belongs either to $\cL_0$ or to $\cL_2$, according the opposite regulus of $\cR$ has $0$ or $2$ lines of $\cW(3, q)$. In this case the opposite regulus of $\cR$ has $0$ or $2$ lines of $\cW(3, q)$ according as the polar of $\langle C \rangle$ with respect to the orthogonal polarity of $\cQ(4, q)$ is a line external or secant to $\cQ(4, q)$. Therefore, in $\cW(3, q)$, $\ell$ can be chosen in $q(q-1)^2/2$ ways such that it belongs to $\cL_0$ and in $(q^3-q)/2$ ways such that it belongs to $\cL_2$. Since there are $q^2(q^2+q+1)$ lines not of $\cW(5, q)$ meeting both $\Pi_1$, $\Pi_2$ in one point, we get $|\cL_0| = q^3(q-1)(q^3-1)/2$ and $|\cL_2| = q^3(q+1)(q^3-1)/2$.    
\end{proof}

\begin{lemma}
Let $\Pi_3$ be a plane of $\cW(5, q)$ skew to $\Pi_1$ and $\Pi_2$. The $q+1$ planes of the Segre variety $\Sigma_{1, 2}$ of $\PG(5, q)$ determined by $\Pi_1$, $\Pi_2$, $\Pi_3$ are generators of $\cW(5, q)$. 
\end{lemma}
\begin{proof}
Let $A$ be an invertible matrix of $\cS_{3, q}$ and consider the symplectic Segre variety $\Sigma_{1, 2}$ determined by the planes $\Pi_1$, $\Pi_2$ and $L(A)$. Direct computations show that the remaining $q-2$ planes of $\Sigma_{1, 2}$ are the planes $L(\lambda A)$, where $\lambda \in \GF(q) \setminus \{0, 1\}$. 
\end{proof}

We will refer to a Segre variety $\Sigma_{1,2}$ of $\PG(5, q)$ whose $q+1$ planes are generators of $\cW(5, q)$ as a {\em symplectic Segre variety} of $\cW(5, q)$. As a consequence of Lemma \ref{segre_symp}, the following corollary arises.

\begin{cor}\label{corollary}
Let $\Sigma_{1, 2}$ be a symplectic Segre variety containing $\Pi_1$, $\Pi_2$ and let $\Pi_3$ be a plane of $\Sigma_{1, 2}$, $\Pi_3 \ne \Pi_1$, $\Pi_3 \ne \Pi_2$. Then $\Pi_3$ contains one line of $\cL_0$, $q+1$ lines of $\cL_1$ and $q^2-1$ lines of $\cL_2$, if $q$ is even, and $q(q-1)/2$ lines of $\cL_0$, $q+1$ lines of $\cL_1$ and $q(q+1)/2$ lines of $\cL_2$, if $q$ is odd.
\end{cor}
\begin{proof}
From Lemma \ref{segre_symp}, there is a non--degenerate polarity $\rho$ of $\Pi_3$. Note that if $\ell$ is a line of $\Pi_3$, then $\ell^\rho = T_\ell^\perp\cap\Pi_3$. If $q$ is even, $\rho$ is a pseudo--polarity and the unique line of $\Pi_3$ belonging to $\cL_0$ is the line $\ell_0$ consisting of its absolute points. The other lines of $\Pi_3$ belong to $\cL_1$ if they pass through $\ell_0^\rho$ and to $\cL_2$ otherwise. 
If $q$ is odd, $\rho$ is an orthogonal polarity and its absolute points form a conic, say $C$. A line of $\Pi_3$ belongs either to $\cL_1$, or to $\cL_0$ or to $\cL_2$ according as it is tangent, external or secant to $\cC$, respectively.
\end{proof}

For a point $P \in \Pi_2$, let $\Sigma_P$ denote a 3--space contained in $P^\perp$ and not containing $P$. When restricted to $\Sigma_P$, the polarity $\perp$ defines a non--degenerate symplectic polar space of $\Sigma_P$, say $\cW_P$. Moreover $r_P = \Sigma_P \cap \Pi_1$ and $t_P = \Sigma_P \cap \Pi_2$ are lines of $\cW_P$. In what follows we investigate the action of the group $G_{\Pi_2}$ on $\cG$. 

\begin{lemma}\label{orbits}
The group $G_{\Pi_2}$ has the following orbits on $\cG$:
\begin{itemize}
\item the plane $\Pi_2$;
\item $\cG_1$ of size $q^3-1$ consisting of the planes of $\cG$ meeting $\Pi_2$ in a line,
\item $\cG_2$ of size $q^3-1$ consisting of the planes of $\cG$ meeting $\Pi_2$ in a point and no plane of $\cG_1$ in a line,
\item $\cG_3$ of size $(q^2-1)(q^3-1)$ consisting of the planes of $\cG$ meeting $\Pi_2$ in a point and $q$ planes of $\cG_1$ in a line,
\item $\cG_4$ of size $q^2(q^3-1)(q-1)$ consisting of the planes of $\cG$ disjoint from $\Pi_2$,
\end{itemize}
if $q$ is even and
\begin{itemize}
\item the plane $\Pi_2$;
\item $\cG_1$ of size $(q^3-1)/2$ consisting of the planes of $\cG$ meeting $\Pi_2$ in a line and having $q^2$ points of $\cP_1$,
\item $\cG_2$ of size $(q^3-1)/2$ consisting of the planes of $\cG$ meeting $\Pi_2$ in a line and having $q^2$ points of $\cP_2$,
\item $\cG_3$ of size $q(q-1)(q^3-1)/2$ consisting of the planes of $\cG$ meeting $\Pi_2$ in a point and $q+1$ planes of $\cG_1 \cup \cG_2$ in a line,
\item $\cG_4$ of size $q(q+1)(q^3-1)/2$ consisting of the planes of $\cG$ meeting $\Pi_2$ in a point and $q-1$ planes of $\cG_1 \cup \cG_2$ in a line,
\item two orbits, say $\cG_5$ and $\cG_6$, both of size $q^2(q^3-1)(q-1)/2$, consisting of planes of $\cG$ disjoint from $\Pi_2$,
\end{itemize}
if $q$ is odd.
\end{lemma}
\begin{proof}
There are $q^3-1$ members of $\cG$ intersecting $\Pi_2$ in a line. The number of generators of $\cW(5,q)$ through $P$ disjoint from $\Pi_1$ and intersecting $\Pi_2$ exactly in $P$ equals the number of lines of $\cW_P$ disjoint from both $r_P$ and $t_P$, and they are $q^2(q-1)$. As the point $P$ varies on $\Pi_2$ we get $q^2(q^3-1)$ generators of $\cW(5,q)$ disjoint from $\Pi_1$ and intersecting $\Pi_2$ at exactly one point. Hence there are $q^2(q-1)(q^3-1)$ generators of $\cW(5,q)$ disjoint from both $\Pi_1$ and $\Pi_2$. Alternatively, if $A \in \cS_{3, q}$, then $L(A) \cap \Pi_2$ is a $(3 - \rk(A) - 1)$--space of $\Pi_1$. 

Let $\pi$ be a plane of $\cG$ intersecting $\Pi_2$ in a line. Since $G_{\Pi_2}$ is transitive on lines of $\Pi_2$, we may assume without loss of generality that $\pi \cap \Pi_2$ is the line $U_2 U_3$. Then $\pi: X_5 = X_6 = 0, X_4 = z X_1$, for some $z \in \GF(q) \setminus \{0\}$. Let $g$ be the projectivity of $G_{\Pi_2}$ associated with the matrix \eqref{matrice_s1}, for some $T \in \GL(3 , q)$. Then $g$ stabilizes $\pi$ if and only if 
$$
T = \begin{pmatrix} 
\pm 1 & * & * \\
0 &  T' & \\
0 &  & \\
\end{pmatrix},
$$
where $T' \in \GL(2, q)$. Hence $|Stab_{G_{\Pi_2}}(\pi)| = gcd(2, q-1) q^2 |\GL(2,q)|/gcd(2, q-1)$ and $|\pi^{G_{\Pi_2}}|$ equals $q^3-1$ if $q$ is even or $(q^3-1)/2$ if $q$ is odd. In the even characteristic case the planes of $\cG$ meeting $\Pi_2$ in a line are permuted in a unique orbit, say $\cG_1$. In the odd characteristic case, there are two orbits, say $\cG_1$ and $\cG_2$; it can be seen that representatives for $\cG_1$ and $\cG_2$ are $\pi_1: X_5 = X_6 = X_4 - z_1 X_1 = 0$ and $\pi_2: X_5 = X_6 = X_4 - z_2 X_1 = 0$, respectively, where $z_1$ is a non--zero square in $\GF(q)$ and $z_2$ is a non--square of $\GF(q)$. Moreover $Stab_{G_{\Pi_2}}(\pi)$ acts transitively on the $q^2$ points of $\pi \setminus \Pi_2$.

Let $\pi$ be a plane of $\cG$ and intersecting $\Pi_2$ in the point $P$ and let $\ell$ be the line of $\cW_P$ obtained by intersecting $\pi$ with $\Sigma_P$. Let $\cR$ be the regulus determined by $r_P, t_P, \ell$ and $\cR^{o}$ be its opposite regulus. 

Assume that $q$ is even. There are $q-1$ possibilities for the line $\ell$ such that the regulus $\cR^{o}$ contains $q+1$ lines of $\cW_P$, and $(q^2-1)(q-1)$ possibilities for $\ell$ such that $\cR^{o}$ contains exactly one line of $\cW_P$. Varying $P$ in $\Pi_2$ we get two sets, namely $\cG_2$ and $\cG_3$, of size $q^3-1$ and $(q^2-1)(q^3-1)$, respectively. Observe that there are $0$ or $q$ planes of $\cG_1$ meeting $\pi$ in a line, according as $\pi$ belongs to $\cG_2$ or $\cG_3$. We claim that $\cG_2$ and $\cG_3$ are two $G_{\Pi_2}$--orbits. Let $\pi$ be the plane with equations $X_2 + X_6 = X_3 + X_5 = X_4 = 0$. Direct computations show that $\pi \in \cG_2$. The projectivity $g$ of $G_{\Pi_2}$ associated with the matrix \eqref{matrice_s1}, $T \in \GL(3 , q)$, stabilizes $\pi$ if and only if 
$$
T = \begin{pmatrix} 
x & 0 & 0 \\
* &  T' & \\
* &  & \\
\end{pmatrix},
$$
where $x \in \GF(q) \setminus \{0\}$ and $T' \in \SL(2, q)$. Hence $|Stab_{G_{\Pi_2}}(\pi)| = q^2 (q-1) |\SL(2,q)|$ and $|\pi^{G_{\Pi_2}}| = q^3-1 = |\cG_2|$. Let $\pi$ be the plane having equations $X_2 + X_5 = X_3 + X_6 = X_4 = 0$. Direct computations show that $\pi \in \cG_3$. The projectivity $g$ of $G_{\Pi_2}$ associated with the matrix \eqref{matrice_s1}, $T \in \GL(3 , q)$, stabilizes $\pi$ if and only if 
$$
T = \begin{pmatrix} 
x & 0 & 0 \\
* & y & y' \\
* & y' & y \\
\end{pmatrix},
$$
where $x, y, y' \in \GF(q)$, $x \ne 0$ and $y^2+y'^2 = 1$. Hence $|Stab_{G_{\Pi_2}}(\pi)| = q^3 (q-1)$ and $|\pi^{G_{\Pi_2}}| = (q^2-1)(q^3-1) = |\cG_3|$.

Assume that $q$ is odd. There are $q(q-1)^2/2$ possibilities for the line $\ell$ such that the regulus $\cR^{o}$ contains no line of $\cW_P$, and $q(q^2-1)/2$ possibilities for $\ell$ such that $\cR^{o}$ contains exactly two lines of $\cW_P$. Varying $P$ in $\Pi_2$ we get two sets, say $\cG_3$ and $\cG_4$ of size $q(q-1)(q^3-1)/2$ and $q(q+1)(q^3-1)/2$ respectively. Observe that there are $q+1$ or $q-1$ planes of $\cG_1 \cup \cG_2$ meeting $\pi$ in a line, according as $\pi$ belongs to $\cG_3$ or $\cG_4$. Again we want to show that $\cG_3$ and $\cG_4$ are two $G_{\Pi_2}$--orbits. Let $\alpha$ be a fixed non--square in $\GF(q)$ and let $\pi$ be the plane with equations $X_5 - \alpha^2 X_2 = X_6 + \alpha X_3 = X_4 = 0$. Direct computations show that $\pi \in \cG_3$. The projectivity $g$ of $G_{\Pi_2}$ associated with the matrix \eqref{matrice_s1}, $T \in \GL(3 , q)$, stabilizes $\pi$ if and only if 
$$
T = \begin{pmatrix} 
x & 0 & 0 \\
* & y & - \alpha y' \\
* & y' & - y \\
\end{pmatrix}
\mbox{ or }
T = \begin{pmatrix} 
x & 0 & 0 \\
* & y & \alpha y' \\
* & y' & y \\
\end{pmatrix}, 
$$
where $x, y, y' \in \GF(q)$, $x \ne 0$ and $y^2 - \alpha y'^2 = 1$. Note that there are $q+1$ couple $(y, y') \in \GF(q) \times \GF(q)$ such that $y^2 - \alpha y'^2 = 1$. Therefore $|Stab_{G_{\Pi_2}}(\pi)| = q^2 (q-1) (q+1)$ and $|\pi^{G_{\Pi_2}}| = q(q-1)(q^3-1)/2 = |\cG_3|$. Let $\pi$ be the plane having equations $X_2 - X_6 = X_3 - X_5 = X_4 = 0$. Direct computations show that $\pi \in \cG_4$. The projectivity $g$ of $G_{\Pi_2}$ associated with the matrix \eqref{matrice_s1}, $T \in \GL(3 , q)$, stabilizes $\pi$ if and only if 
$$
T = \begin{pmatrix} 
x & 0 & 0 \\
* & y & 0 \\
* & 0 & y^{-1} \\
\end{pmatrix}
\mbox{ or }
T = \begin{pmatrix} 
x & 0 & 0 \\
* & 0 & y \\
* & y^{-1} & 0 \\
\end{pmatrix}, 
$$
where $x, y \in \GF(q) \setminus \{0\}$. In this case $|Stab_{G_{\Pi_2}}(\pi)| = q^2 (q-1)^2$ and $|\pi^{G_{\Pi_2}}| = q(q+1)(q^3-1)/2 = |\cG_4|$.

From Lemma \ref{dis}, the group $G_{\Pi_2}$ has one or two orbits on generators of $\cW(5, q)$ skew to $\Pi_1$ and $\Pi_2$. 
\end{proof}

\begin{lemma}\label{lemma1}
Assume that $q$ is even. Let $\Pi \in \cG_2$ and $\Pi' \in \cG_3$. Then the number of planes of $\cG_2$ meeting $\Pi$ or $\Pi'$ in a line is zero or one, whereas the number of planes of $\cG_3$ meeting $\Pi$ or $\Pi'$ in a line equals $q^2-1$ or $q^2-q-2$.
\end{lemma}

\begin{proof}
Let $P = \Pi \cap \Pi_2$, $P' = \Pi' \cap \Pi_2$, $\ell = \Pi \cap \Sigma_P$, $\ell' = \Pi' \cap \Sigma_{P'}$, $\cR$ be the regulus determined by $r_P, t_P, \ell$ and $\cR'$ be the regulus determined by $r_P', t_P', \ell'$. From the proof of Lemma \ref{orbits}, the opposite regulus of $\cR$, say $\cR^o$, consists of lines of $\cW_P$, whereas the opposite regulus of $\cR'$, say $\cR'^o$, has exactly one line of $\cW_P'$. 

If a plane $\gamma \in \cG_2 \cup \cG_3$ intersects $\Pi$ in a line, then $\gamma = \langle P, s \rangle$, where $s$ is a line of $\cW_P$ intersecting $\ell$ and skew to $r_P$ and $t_P$. Moreover, $\gamma$ belongs to $\cG_2$ or $\cG_3$ according as there are $q+1$ or one line of $\cW_P$ meeting $s, r_P, t_P$. Since the number of lines of $\cW_P$ intersecting $\ell$ at a point and skew to $r_P$ and $t_P$ equals $q^2-1$ and, if $s$ is one of these lines, there is exactly one line of $\cW_P$ meeting $s, r_P, t_P$, the statement holds true in this case.

Similarly, if a plane $\gamma \in \cG_2 \cup \cG_3$ intersects $\Pi'$ in a line, then $\gamma = \langle P', s \rangle$, where $s$ is a line of $\cW_{P'}$ intersecting $\ell'$ and skew to $r_P'$ and $t_P'$. Moreover, $\gamma$ belongs to $\cG_2$ or $\cG_3$ according as there are $q+1$ or one line of $\cW_{P'}$ meeting $s, r_P', t_P'$. The point line dual of $\cW_{P'}$ is a parabolic quadric $\cQ(4, q)$, the lines $r_P'$, $t_P'$ and $\ell'$ correspond to three points, say $R, T, L$, such that the line $R T$ meets $\cQ(4, q)$ only in $R$ and $T$. Moreover, the regulus $\cR'$ corresponds to a conic $C$ of $\cQ(4, q)$, where $R, T, L \in C$. Let $N$ be the nucleus of $\cQ(4, q)$. Then $N$ does not belong to the plane $\langle C \rangle$ and the points $R, T, N$ span a plane meeting $\cQ(4, q)$ in a conic, say $C'$. Hence there is a plane $\gamma \in \cG_2$ meeting $\Pi'$ in a line if and only if there is a point $U$ of $C'$ such that the line $U L$ is a line of $\cQ(4, q)$. There exists only one such a point: the intersection point between the three--space containing the lines of $\cQ(4, q)$ through $L$ and the conic $C'$. Analogously, there is a plane $\gamma \in \cG_3$ meeting $\Pi'$ in a line if and only if there is a point $U \in \cQ(4, q)$ not belonging to $C'$ such that the line $U L$ is a line of $\cQ(4, q)$ and the lines $U R$ and $U T$ are not lines of $\cQ(4, q)$. There exist exactly $q^2-q-2$ points having these properties.
\end{proof}

\begin{lemma}\label{lemma2}
Assume that $q$ is odd. Let $\Pi \in \cG_5$ and $\Pi' \in \cG_6$, then either $|\Pi \cap \cP_1| = |\Pi' \cap \cP_2| = q(q-1)/2$ and $|\Pi \cap \cP_2| = |\Pi' \cap \cP_1| = q(q+1)/2$ or $|\Pi \cap \cP_1| = |\Pi' \cap \cP_2| = q(q+1)/2$ and $|\Pi \cap \cP_2| = |\Pi' \cap \cP_1| = q(q-1)/2$. Moreover through a line of $\cL_0 \cup \cL_2$, there pass $(q-1)/2$ planes of $\cG_5$ and $(q-1)/2$ planes of $\cG_6$, whereas the $q$ generators passing through a line of $\cL_1$ and skew to $\Pi_1$, $\Pi_2$ are planes either of $\cG_5$ or of $\cG_6$. 
\end{lemma}
\begin{proof}
Let $A$ be an invertible matrix of $\cS_{3, q}$ and consider the symplectic Segre variety $\Sigma_{1, 2}$ determined by the $\Pi_1$, $\Pi_2$ and $L(A)$. The planes of $\Sigma_{1, 2}$ distinct from $\Pi_1$ and $\Pi_2$ are $L(\lambda A)$, where $\lambda \in \GF(q) \setminus \{0\}$. Taking into account Remark \ref{oss1}, we have that $(q-1)/2$ members of $\Sigma_{1, 2}$ belong to $\cG_{5}$ and $(q-1)/2$ members of $\Sigma_{1, 2}$ belong to $\cG_{6}$. Moreover, taking into account Remark \ref{oss}, if $\lambda$ is a non--zero square of $\GF(q)$, then the point of $L(A)$ given by $(x, y, z, 0, 0, 0) + (0, 0, 0, x, y, z) \begin{pmatrix} 0_3 & A \end{pmatrix}^t$ belongs to $\cP_1$ if and only if the point of $L(\lambda A)$ given by $(x, y, z, 0, 0, 0) + (0, 0, 0, x, y, z) \lambda \begin{pmatrix} 0_3 & A \end{pmatrix}^t$ belongs to $\cP_1$, whereas if $\lambda$ is a non--square of $\GF(q)$, then the point of $L(A)$ given by $(x, y, z, 0, 0, 0) + (0, 0, 0, x, y, z) \begin{pmatrix} 0_3 & A \end{pmatrix}^t$ belongs to $\cP_1$ if and only if the point of $L(\lambda A)$ given by $(x, y, z, 0, 0, 0) + (0, 0, 0, x, y, z) \lambda \begin{pmatrix} 0_3 & A \end{pmatrix}^t$ belongs to $\cP_2$. Let $\Pi = L(\lambda A) \in \cG_5$, $\Pi' = L(\lambda' A) \in \cG_6$. From Lemma \ref{segre_symp}, there is a non--degenerate conic $C$ (resp. $C'$) of $\Pi$ (resp. $\Pi'$). Observe that exactly one of the two following possibilities occurs: either $\Pi \cap \cP_1$ are the points of $\Pi$ internal to $C$, $\Pi \cap \cP_2$ are the points of $\Pi$ external to $C$, $\Pi' \cap \cP_2$ are the points of $\Pi'$ internal to $C'$, $\Pi' \cap \cP_1$ are the points of $\Pi'$ external to $C'$, or $\Pi \cap \cP_1$ are the points of $\Pi$ external to $C$, $\Pi \cap \cP_2$ are the points of $\Pi$ internal to $C$, $\Pi' \cap \cP_2$ are the points of $\Pi'$ external to $C'$, $\Pi' \cap \cP_1$ are the points of $\Pi'$ internal to $C'$. 

Let $\ell$ be a line of $\cW(5, q)$ and let $\Pi_3$ be a generator of $\cW(5, q)$ skew to $\Pi_1$, $\Pi_2$ such that $\ell \subset \Pi_3$. Denote by $\rho$ the non--degenerate polarity of $\Pi_3$ arising from Lemma \ref{segre_symp} and let $C_3$ be the corresponding non--degenerate conic. If $\ell \in \cL_0 \cup \cL_2$, then $T_{\ell}^\perp$ is a line meeting both $\Pi_1$, $\Pi_2$ in a point and it is not a line of $\cW(5, q)$. Hence, by Remark \ref{oss}, $T_{\ell}^\perp$ contains $(q-1)/2$ points of $\cP_1$ and $(q-1)/2$ points of $\cP_2$. Since $\ell^\rho$ is the point $\Pi_3 \cap T_\ell^\perp$, we have that through $\ell$ there pass $(q-1)/2$ planes of $\cG_5$ and $(q-1)/2$ planes of $\cG_6$. If $\ell \in \cL_1$, then $T_{\ell}^\perp$ is a line of $\cW(5, q)$ meeting both $\Pi_1$, $\Pi_2$ in a point. In this case $\ell \cap C_3$ consists of one point, say $Q$. The $q$ points of $\ell$ distinct from $Q$ are external to $C_3$ and hence they all lie in a unique point--orbit, that is either $\cP_1$ or $\cP_2$. Therefore the $q$ generators of $\cW(5, q)$ passing through $\ell$ and skew to $\Pi_1$, $\Pi_2$ are such that the external points of their corresponding conics are all points belonging to the same point orbit, that is either $\cP_1$ or $\cP_2$. Hence these $q$ generators lie in the same $G_{\Pi_2}$--orbit which contains $\Pi_3$.   
\end{proof}

\section{Hermitian matrices and Hermitian polar spaces}

Let $\omega \in \GF(q^2) \setminus \{0\}$ such that $\omega^q = -\omega$ and let $\cH(2n-1, q^2)$ be the non--degenerate Hermitian polar space of $\PG(2n-1, q^2)$ associated with the following Hermitian form 
$$
\begin{pmatrix} X_1, \dots, X_{2n} \end{pmatrix}
\begin{pmatrix}
0_n & \omega I_n\\
\omega^q I_n & 0_n
\end{pmatrix}
\begin{pmatrix} Y_1^q \\ \vdots \\ Y_{2n}^q \end{pmatrix}.
$$ 
Let $\PGU(2n, q^2)$ be the group of projectivities of $\PG(2n-1, q^2)$ stabilizing $\cH(2n-1, q^2)$. Denote by $\Lambda_1$ the $(n-1)$--space of $\PG(2n -1, q^2)$ spanned by $U_{n+1}, \dots, U_{2n}$. Then $\Lambda_1$ is a generator of $\cH(2n-1, q^2)$. Denote by $\bar{G}$ the stabilizer of $\Lambda_1$ in $\PGU(2n, q^2)$. In this case an element of $\bar{G}$ is represented by the matrix
\begin{equation}
\begin{pmatrix} \label{matrice_h}
T^{-t} & 0_n \\
H_0^q T^{-t} & T^q \\
\end{pmatrix},
\end{equation}
where $T \in \GL(n, q^2)$ and $H_0 \in \cH_{n, q^2}$. Hence $\bar{G} \simeq \cH_{n, q^2} \rtimes \left(\GL(n, q^2)/\langle a I_{2n} \rangle\right)$, where $a^{q+1} = 1$, and  
$$
|\bar{G}| = \frac{q^{n^2} \prod_{i = 0}^{n-1} (q^{2n}-q^{2i})}{q+1}.
$$
Define an action of $\cH_{n, q^2} \rtimes \GL(n, q^2)$ on $\cH_{n, q^2}$ as follows
$$
\left((H_0, T), H \right) \in \left(\cH_{n, q^2} \rtimes \GL(n, q^2)\right) \times \cH_{n, q^2} \longmapsto T H \left(T^q\right)^t + H_0 \in \cH_{n, q^2}.
$$
Its orbitals are the relations of an association scheme, the so called association scheme of Hermitian matrices \cite{W, S2}. As in the symmetric case there is a correspondence between $\cH_{n, q^2}$ and the set $\bar{\cG}$ of generators of $\cH(2n-1, q^2)$ disjoint from $\Lambda_1$ (see also \cite[Proposition 9.5.10]{BCN}). The proof is similar to that of the symmetric case and hence it is omitted.
\begin{lemma}\label{Hermitian}
There is a bijection between $\cH_{n, q^2}$ and $\bar{\cG}$ such that $\cH_{n, q^2} \rtimes \GL(n, q^2)$ acts on $\cH_{n, q^2}$ as $\bar{G}$ acts on $\bar{\cG}$. In particular, a $d$--code of $\cH_{n, q^2}$ corresponds to a set of generators of $\cH(2n-1, q^2)$ disjoint from $\Lambda_1$ pairwise intersecting in at most an $(n - d - 1)$--space, and conversely.
\end{lemma}
As before, if $\Lambda_2 = L(0_n)$, the previous lemma implies that the number of orbits of $\bar{G}_{\Lambda_2}$ on $\bar{\cG}$ equals the number of relations of the association scheme on Hermitian matrices.

\begin{lemma}[\cite{Thas1}]\label{segre_herm}
Let $\Pi_3$ be a generator of $\cH(2n-1, q^2)$ disjoint from $\Pi_1$ and $\Pi_2$. The points $P \in \Pi_3$ such that there exists a line of $\cH(2n-1, q^2)$ through $P$ intersecting $\Pi_1$ and $\Pi_2$ are the absolute points of a non--degenerate unitary polarity of $\Pi_3$. 
\end{lemma}

\section{$2$--codes of $\cS_{n, q}$ or $\cH_{n, q^2}$}

Let $\Gamma_{\cW}$ or $\Gamma_{\cH}$ be the graph whose vertices are the generators of $\cW(2n-1, q)$ or $\cH(2n-1, q^2)$ and two vertices are adjacent whenever they meet in an $(n-2)$--space. Then $\Gamma_{\cW}$ or $\Gamma_{\cH}$ is a distance regular graph having diameter $n$, see \cite[Section 9.4]{BCN}. A coclique of $\Gamma_{\cW}$ or $\Gamma_{\cH}$ is a set of generators of $\cW(2n-1, q)$ or $\cH(2n-1, q^2)$ pairwise intersecting in at most an $(n-3)$--space. 
\begin{lemma}[Theorem 9.4.3, \cite{BCN}]
The eigenvalues $\theta_j$, $0 \le j \le n$, of $\Gamma_{\cW}$ are:
$$
\theta_j = \frac{q^j(q^{n-2j+1}-1)}{q-1}-1,  \mbox{ with multiplicity\ \ }  f_j = \frac{q^j(q^{n-2j+1} + 1)}{q^{n-j+1} + 1} \prod_{i=1}^{j} \frac{q^{2(n-i+1)} - 1}{(q^i - 1)(q^{i-1} + 1)} .
$$
The eigenvalues $\lambda_j$, $0 \le j \le n$, of $\Gamma_{\cH}$ are:
$$
\lambda_j = \frac{q^{2j}(q^{2n-4j+1}-1)}{q^2-1}-\frac{1}{q+1},  \mbox{ with multiplicity\ \ }  g_j = \frac{q^{2j}(q^{2n-4j+1} + 1)}{q^{2n-2j+1} + 1} \prod_{i=1}^{j} \frac{(q^{2n-2i+2} - 1)(q^{2n-2i+1} + 1)}{(q^{2i} - 1)(q^{2i-1} + 1)} .
$$
\end{lemma}
The eigenvalue $\theta_j$ is positive or negative, according as $0 \le j \le \left\lceil\frac{n-1}{2}\right\rceil$ or $\left\lceil\frac{n-1}{2}\right\rceil + 1 \le j \le n$, respectively, and
$$
\deg(f_j) = 
\begin{cases} 
n j + j (n - 2j +1) & \mbox{ if  } 0 \le j \le \left\lceil\frac{n-1}{2}\right\rceil, \\
n j + (j - 1) (n - 2j +1) & \mbox{ if  } \left\lceil\frac{n-1}{2}\right\rceil + 1 \le j \le n. \\
\end{cases}
$$
Moreover $\deg(f_i) < \deg(f_j)$, if $0 \le i < j \le \left\lceil\frac{n-1}{2}\right\rceil$, and $\deg(f_i) > \deg(f_j)$, if $\left\lceil\frac{n-1}{2}\right\rceil + 1 \le i < j \le n$. From the Cvetkovi{\'c} bound (Lemma \ref{cvetkovic}), it follows that $\alpha(\Gamma_{\cW}) \le \sum_{j=0}^{\frac{n-1}{2}} f_j$, if $n$ is odd. Note that the Hoffman bound gives a better upper bound for $\alpha(\Gamma_{\cW})$ than the Cvetkovi{\'c} bound in the case when $n$ is even. Hence, the following result arises.     

\begin{theorem}\label{upper_symp}
Let $\cX$ be a set of generators of $\cW(2n-1, q)$ pairwise intersecting in at most an $(n-3)$--space. Then
$$
|\cX| \le 
\begin{cases}
\sum_{j = 0}^{\frac{n-1}{2}} \frac{q^j(q^{n-2j+1} + 1)}{q^{n-j+1} + 1} \prod_{i=1}^{j} \frac{q^{2(n-i+1)} - 1}{(q^i - 1)(q^{i-1} + 1)} & \mbox{ if } n \mbox{ is odd, } \\ 
\prod_{i = 2}^{n} (q^i + 1) & \mbox{ if } n \mbox{ is even. }
\end{cases}
$$
\end{theorem}

From Lemma \ref{symplectic}, the size of the largest $2$--codes of $\cS_{n, q}$ coincides with the maximum number of generators of $\cW(2n-1, q)$ disjoint from $\Pi_1$ such that they pairwise intersect in at most an $(n-3)$--space. Hence we have the following corollary.

\begin{cor}\label{upper_bound_sym}
Let $\cC$ be a $2$--code of $\cS_{n, q}$, then 
$$
|\cC| \le 
\begin{cases}
\sum_{j = 0}^{\frac{n-1}{2}} \frac{q^j(q^{n-2j+1} + 1)}{q^{n-j+1} + 1} \prod_{i=1}^{j} \frac{q^{2(n-i+1)} - 1}{(q^i - 1)(q^{i-1} + 1)} & \mbox{ if } n \mbox{ is odd, } \\ 
\prod_{i = 2}^{n} (q^i + 1) & \mbox{ if } n \mbox{ is even. }
\end{cases}
$$
\end{cor}
The term of highest degree in Corollary \ref{upper_bound_sym} is $q^{n(n+1)/2-1}/2$ if $n$ is odd and $q^{n(n+1)/2-1}$ if $n$ is even. The previous known upper bound for the size of a $2$--code of $\cS_{n, q}$ was $q^{n(n-1)/2+1}(q^{n-1}+1)/(q+1)$ for $q$ odd \cite[Proposition 3.7]{S1}, and $q^{n(n+1)/2} - q^n + 1$ for $q$ even \cite[Proposition 3.4]{S}. Therefore Corollary~\ref{upper_bound_sym} provides better upper bounds if $n$ is odd or if $n$ and $q$ are both even. Regarding $2$--codes of $\cS_{3, q}$, a further improvement will be obtained in Section \ref{case3}.

\bigskip

In the Hermitian case, $\lambda_j$ is positive or negative, according as $0 \le j \le \left\lceil\frac{n-1}{2}\right\rceil$ or $\left\lceil\frac{n-1}{2}\right\rceil + 1 \le j \le n$, respectively, and
$$
\deg(g_j) = 
\begin{cases} 
4 j (n - j) & \mbox{ if  } 0 \le j \le \left\lceil\frac{n-1}{2}\right\rceil, \\
4 j (n - j) - (2n - 4j +1) & \mbox{ if  } \left\lceil\frac{n-1}{2}\right\rceil + 1 \le j \le n. \\
\end{cases}
$$
Moreover $\deg(g_i) < \deg(g_j)$, if $0 \le i < j \le \left\lceil\frac{n-1}{2}\right\rceil$, and $\deg(g_i) > \deg(g_j)$, if $\left\lceil\frac{n-1}{2}\right\rceil + 1 \le i < j \le n$. From the Cvetkovi{\'c} bound, it follows that $\alpha(\Gamma_{\cH}) \le \sum_{j=0}^{\frac{n-1}{2}} g_j$, if $n$ is odd and $\alpha(\Gamma_{\cH}) \le \sum_{\frac{n}{2}+1}^{n} g_j$, if $n$ is even.

\begin{theorem}\label{upper_her}
Let $\cX$ be a set of generators of $\cH(2n-1, q^2)$ pairwise intersecting in at most an $(n-3)$--space. Then 
$$
|\cX| \le 
\begin{cases}
\sum_{j = 0}^{\frac{n-1}{2}}  \frac{q^{2j}(q^{2n-4j+1} + 1)}{q^{2n-2j+1} + 1} \prod_{i=1}^{j} \frac{(q^{2n-2i+2} - 1)(q^{2n-2i+1} + 1)}{(q^{2i} - 1)(q^{2i-1} + 1)} & \mbox{ if } n \mbox{ is odd, }\\
\sum_{\frac{n}{2}+1}^{n} \frac{q^{2j}(q^{2n-4j+1} + 1)}{q^{2n-2j+1} + 1} \prod_{i=1}^{j} \frac{(q^{2n-2i+2} - 1)(q^{2n-2i+1} + 1)}{(q^{2i} - 1)(q^{2i-1} + 1)} & \mbox{ if } n \mbox{ is even. }\\
\end{cases}
$$
\end{theorem}

From Lemma \ref{Hermitian}, the size of the largest $2$--codes of $\cH_{n, q^2}$ coincides with the maximum number of members of $\bar{\cG}$, the set of generators of $\cH(2n-1, q^2)$ disjoint from $\Lambda_1$, such that they pairwise intersect in at most an $(n-3)$--space. Let $\Gamma_{\cH}'$ be the induced subgraph of $\Gamma_{\cH}$ on $\bar{\cG}$. The graph $\Gamma_{\cW}'$ is also known as the last subconstituent of $\Gamma_{\cW}$ or the {\em Hermitian forms graph}, see \cite[Proposition 9.5.10]{BCN}. In particular $\Gamma_{\cW}'$ is a distance--regular graph of diameter $n$. The eigenvalues  of $\Gamma_{\cH}'$ are \cite[Corollary 8.4.4]{BCN}:
$$
\frac{(-q)^{2n-j}-1}{q+1} ,  \mbox{ with multiplicity\ \ }   \prod_{i=1}^{j} \frac{(-q)^{n+1-i} - 1}{(-q)^i - 1} \prod_{i = 0}^{j-1}\left( -(-q)^{n} - (-q)^i\right),
$$
respectively, with $0 \le j \le n$.

In this case the Hoffman bound gives a better upper bound for $\alpha(\Gamma_{\cH}')$ than the Cvetkovi{\'c} bound, that is $q^{(n-1)^2}(q^{2n-1}+1)/(q+1)$. However this was already known \cite[Theorem 2]{S3}. 

\section{Large non--additive rank distance codes}\label{case3}

In the case when $\cC$ is additive much better bounds can be provided. Indeed in \cite[Theorem 4.3]{S2}, \cite[Theorem 1]{S3} the author proved that the largest additive $d$--codes of $\cS_{n, q}$ have size at most either $q^{n(n-d+2)/2}$ or $q^{(n+1)(n-d+1)/2}$, according as $n-d$ is even or odd, respectively, whereas the size of the largest additive $d$--codes of $\cH_{n, q^2}$ cannot exceed $q^{n(n-d+1)}$. As far as regard codes whose size is larger than the additive bound not much is known. In the Hermitian case, if $n$ is even, there is an $n$--code of $\cH_{n, q^2}$ of size $q^n+1$ \cite[Theorem 6]{S3}, \cite[Theorem 18]{GLSV}. Observe that from Lemma \ref{symplectic} and Lemma \ref{Hermitian}, an $n$--code $\cC$ of $\cS_{n, q}$ or $\cH_{n, q^2}$ exists if and only if there exists a partial spread of $\cW(2n-1, q)$ or $\cH(2n-1, q^2)$ of size $|\cC|+1$. It is well--known that the points of $\cW(2n -1, q)$ can be partitioned into $q^n+1$ pairwise disjoint generators of $\cW(2n-1, q)$, that is, $\cW(2n-1, q)$ admits a spread. On the other hand, $\cH(2n-1, q^2)$ has no spread. If $n$ is odd, an upper bound for the largest partial spreads of $\cH(2n-1, q^2)$ is $q^n+1$ \cite{V} and there are examples of partial spreads of that size \cite{L}. If $n$ is even the situation is less clear: upper bounds can be found in \cite{ISX}, as for lower bounds there is a partial spread of $\cH(2n-1, q^2)$ of size $(3q^2-q)/2 + 1$ for $n = 2$, $q > 13$, \cite[p. 32]{ACE} and of size $q^n+2$ for $n \ge 4$ \cite{GLSV}. Here, generalizing the partial spread of $\cH(3, q^2)$, we show the existence of a partial spread of $\cH(2n-1, q^2)$, in the case when $n$ is even and $n/2$ is odd, of size $\left(3q^{n}-q^{n/2}\right)/2 + 1$, see Theorem \ref{partial_spread}. Hence the following result holds true. 

\begin{theorem}
If $n$ is even and $n/2$ is odd, then there exists an $n$--code of $\cH_{n, q^2}$ of size $\frac{3q^{n}-q^{n/2}}{2}$.
\end{theorem}

For small values of $d$, $q$ and $n$, in \cite{S} there are several $d$--codes of $\cS_{n, q}$ and $\cH_{n, q^2}$ whose sizes are larger than the corresponding additive bounds, namely a $2$--code of $\cS_{3, 2}$ of size $22$, a $2$--code of $\cS_{3, 3}$ of size $135$, a $2$--code of $\cS_{3, 4}$ of size $428$, a $2$--code of $\cS_{3, 5}$ of size $934$, a $2$--code of $\cS_{3, 7}$ of size $3100$, a $2$--code of $\cS_{4, 2}$ of size $320$, a $4$--code of $\cS_{5, 2}$ of size $96$, a $2$--code of $\cH_{3, 4}$ of size $120$ and a $4$--code of $\cH_{4, 4}$ of size $37$. Besides these few examples, no $d$--codes whose sizes are larger than the largest possible additive $d$--codes are known.  

\bigskip

In the remaining part of this section we focus on the case $n = 3$. From Corollary \ref{upper_bound_sym}, a $2$--code $\cC$ of $\cS_{3, q}$ has size at most 
$$
\frac{q(q^2+1)(q^2+q+1)}{2} + 1.
$$
First, we improve on the upper bound of the size of $\cC$. Then we construct $2$--codes of $\cS_{3, q}$ and $\cH_{3, q^2}$ that are larger than the largest possible additive $2$--codes. This provides an answer to a question posed in \cite[Section 7]{Sheekey}, see also \cite[p. 176]{S2}. The main results are summarized in the following theorem.
\begin{theorem}
Let $\cC$ be a maximum $2$--code of $\cS_{3, q}$, $q >2$, then 
$$
q^4+q^3+1 \le |\cC| \le \frac{q(q^2-1)(q^2+q+1)}{2} + 1.   
$$
Let $\cC$ be a maximum $2$--code of $\cH_{3, q^2}$, then 
$$
q^6+ \frac{q(q-1)(q^4+q^2+1)}{2} \le |\cC| \le q^4(q^4-q^3+q^2-q+1). 
$$
\end{theorem}

\subsection{Partial spread of $\cH(8m - 5, q^2)$}

Let us consider the projective line $\PG(1, q^{4m-2})$ whose underlying vector space is $V(2, q^{4m-2})$, and let $\cH(1, q^{4m-2})$ be a non--degenerate Hermitian polar space of $\PG(1, q^{4m-2})$ associated with $h$, where $h$ is a sesquilinear form on $V(2, q^{4m-2})$. The vector space $V(2, q^{4m-2})$ can be regarded as a $(4m-2)$--dimensional vector space over $\GF(q^2)$, say $\bar{V}$. More precisely 
$$
\bar{V} = \left\{\left(x, x^{q^2}, \dots, x^{q^{4m-4}}, y, y^{q^2}, \dots, y^{q^{4m-4}}\right) \;\; | \;\; (x, y) \in V(2, q^{4m-2})\right\}.
$$
Let $\PG(4m-3, q^2)$ be the projective space whose underlying vector space is $\bar{V}$ and let 
$$
Tr_{q^{4m-2}|q^2} : x \in \GF(q^{4m-2}) \longmapsto \sum_{i=0}^{2m-2} x^{q^{2i}} \in \GF(q^2)
$$ 
denote the usual {\em trace function}. Note that 
$$
\bar{h} = Tr_{q^{4m-2|q^2}} \circ h : \bar{V} \times \bar{V} \longrightarrow \GF(q^2)
$$ 
is a non--degenerate sesquilinear form on $\bar{V}$ and hence there is a non--degenerate polar space $\cH(4m-3, q^2)$ of $\PG(4m-3, q^2)$ associated with $\bar{h}$. See \cite{Gill} for more details. Let $\rho$ be the unitary polarity of $\PG(4m-3, q^2)$ defining $\cH(4m-3, q^2)$.  

\begin{lemma}\label{pre}
There exists a $(2m-2)$--spread $\bS$ of $\PG(4m-3, q^2)$, such that $q^{2m-1}+1$ members of $\bS$ are generators of $\cH(4m-3, q^2)$ and the remaining $q^{4m-2}-q^{2m-1}$ are such that they occur in $\left(q^{4m-2}-q^{2m-1}\right)/2$ pairs of type $\left\{\Delta, \Delta^{\rho}\right\}$, where $|\Delta \cap \Delta^{\rho}| = 0$.
\end{lemma}
\begin{proof}
With the notation introduced above, if $W$ is a vector subspace of $V(2, q^{4m-2})$ of dimension one, then
$$
\left\{\left(x, x^{q^2}, \dots, x^{q^{4m-4}}, y, y^{q^2}, \dots, y^{q^{4m-4}}\right) \;\; | \;\; (x, y) \in W \right\}
$$ 
is a $(2m-1)$--dimensional vector subspace of $\bar{V}$. Hence a point of $\PG(1, q^{4m-2})$ is sent to a $(2m-2)$--space of $\PG(4m-3, q^2)$ and two distinct $(2m-2)$--spaces of $\PG(4m - 3, q^2)$ so obtained are pairwise skew. Let $\bS$ be the set of $(2m-2)$--spaces of $\PG(4m-3, q^2)$ constructed in this way. Then $\bS$ is a $(2m-2)$--spread of $\PG(4m-3, q^2)$ and $|\bS| = q^{4m-2}+1$. 

Note that $\cH(1, q^{4m-2})$ consists of $q^{2m-1}+1$ points. The polarity of $\PG(1, q^{4m-2})$ defining $\cH(1, q^{4m-2})$ fixes each of these $q^{2m-1}+1$ points and interchanges in pairs the remaining $q^{4m-2}-q^{2m-1}$ points of $\PG(1, q^{4m-2})$. Therefore an element $\Delta$ of $\bS$ is a generator of $\cH(4m-3, q^2)$ if $\Delta$ corresponds to a point of $\cH(1, q^{4m-2})$; otherwise $|\Delta \cap \Delta^{\rho}| = 0$.  
\end{proof}

Let $\cH(8m-5, q^2)$ be a non--degenerate Hermitian polar space of $\PG(8m-5, q^2)$ and let $\perp$ be the unitary polarity of $\PG(8m-5, q^2)$ defining $\cH(8m-5, q^2)$.

\begin{theorem}\label{partial_spread}
$\cH(8m-5, q^2)$ has a partial spread of size $\frac{3q^{4m-2}-q^{2m-1}}{2}+1$.
\end{theorem}
\begin{proof}
Let $\Pi_1$, $\Pi_2$, $\Pi_3$ be three pairwise disjoint generators of $\cH(8m-5, q^2)$. Then $\Pi_i \simeq \PG(4m-3, q^2)$. Moreover, there is a non--degenerate unitary polarity $\rho_i$ of $\Pi_i$, see Lemma \ref{segre_herm}. Let $\cH_{i}$ be the Hermitian polar space of $\Pi_i$, defined by $\rho_i$, $i = 1,2,3$. From Lemma \ref{pre}, there exists a $(2m-2)$--spread $\bS$ of $\Pi_1$, such that $q^{2m-1}+1$ members of $\bS$ are generators of $\cH_1$ and the remaining $q^{4m-2}-q^{2m-1}$ are such that they occur in $\left(q^{4m-2}-q^{2m-1}\right)/2$ pairs of type $\left\{\Delta, \Delta^{\rho_1}\right\}$, where $|\Delta \cap \Delta^{\rho_1}| = 0$. For an element $\Delta_1$ of $\bS$, let $\Delta_2 = \langle \Pi_3, \Delta_1 \rangle \cap \Pi_2$ and $\Delta_3 = \langle \Pi_2, \Delta_1 \rangle \cap \Pi_3$. Then $\Delta_2^{\rho_2} = \Delta_1^\perp \cap \Pi_2 = \Delta_3^\perp \cap \Pi_2$ and $\Delta_3^{\rho_3} = \Delta_1^\perp \cap \Pi_3 = \Delta_2^\perp \cap \Pi_3$. Similarly, $\Delta_1^{\rho_1} = \Delta_2^\perp \cap \Pi_1 = \Delta_3^\perp \cap \Pi_1$. 

If $\Delta_1$ is a generator of $\cH_1$, then $\Delta_i^{\rho_i} = \Delta_i$, $i = 1,2,3$, and $\langle \Delta_1, \Delta_2 \rangle$ is a generator of $\cH(8m-5, q^2)$. Varying $\Delta_1$ among the $q^{2m-1}+1$ members of $\bS$ that are generators of $\cH_1$, one obtains a set $Z_1$ of $q^{2m-1}+1$ generators of $\cH(8m-5, q^2)$ that are pairwise disjoint. 

If $\Delta_1$ is such that $|\Delta_1 \cap \Delta_1^{\rho_1}| = 0$, then consider the following generators of $\cH(8m-5, q^2)$:
$$
\left\{ \left\langle \Delta_1, \Delta_2^{\rho_2} \right\rangle, \left\langle \Delta_2^{\rho_2}, \Delta_3 \right\rangle, \left\langle \Delta_3, \Delta_1^{\rho_1} \right\rangle, \left\langle \Delta_1^{\rho_1}, \Delta_2 \right\rangle, \left\langle \Delta_2, \Delta_3^{\rho_3} \right\rangle, \left\langle \Delta_3^{\rho_3}, \Delta_1 \right\rangle \right\}.
$$
Among these six generators, we can always choose three of them such that they are pairwise disjoint. For instance 
$$
\left\langle \Delta_1, \Delta_2^{\rho_2} \right\rangle, \left\langle \Delta_3, \Delta_1^{\rho_1} \right\rangle, \left\langle \Delta_2, \Delta_3^{\rho_3} \right\rangle
$$
are pairwise skew since any two of them span the whole $\PG(8m-5, q^2)$. Repeating this process for each of the $\left(q^{4m-2}-q^{2m-1}\right)/2$ couples $\left\{\Delta_1, \Delta_1^{\rho_1}\right\}$ such that $|\Delta_1 \cap \Delta_1^{\rho_1}| = 0$, a set $Z_2$ of $3\left(q^{4m-2}-q^{2m-1}\right)/2$ pairwise disjoint generators of $\cH(8m-5, q^2)$ is obtained. Again two members of $Z_1 \cup Z_2$ span the whole $\PG(8m-5, q^2)$. Therefore $Z_1 \cup Z_2$ is a partial spread of $\cH(8m-5, q^2)$ of size $3\left(q^{4m-2}-q^{2m-1}\right)/2 + q^{2m-1}+1 = \left(3q^{4m-2}-q^{2m-1}\right)/2+1$.
\end{proof}

\subsection{$2$--codes of $\cS_{3, q}$}

Let $\perp$ be the symplectic polarity of $\PG(5, q)$ defining $\cW(5, q)$. Recall that $\cG$ is the set of $q^6$ planes of $\cW(5, q)$ that are disjoint from $\Pi_1$, the group $G$ is the stabilizer of $\Pi_1$ in $\PSp(6, q)$, $\Pi_2 = L(0_3)$, and $G_{\Pi_2}$ is the stabilizer of $\Pi_2$ in $G$. For a point $P$ in $\Pi_2$, let $\Sigma_P$ denote a $3$--space contained in $P^\perp$ and not containing $P$. When restricted to $\Sigma_P$, the polarity $\perp$ defines a non--degenerate symplectic polar space of $\Sigma_P$, say $\cW_{P}$. Moreover $r_P = \Sigma_P \cap \Pi_1$ and $t_P = \Sigma_P \cap \Pi_2$ are lines of $\cW_P$. 

\subsubsection{The upper bound}

The graph $\Gamma_{\cW}$ has valency $q(q^2+q+1)$. Let $\Gamma_{\cW}'$ be the induced subgraph of $\Gamma_{\cW}$ on $\cG$. The graph $\Gamma_{\cW}'$ is also known as the {\em second subconstituent} of $\Gamma_{\cW}$, see \cite{CK}. Then $\Gamma_{\cW}'$ is connected, has valency $q^3-1$ and it is vertex--transitive, since $G \le Aut(\Gamma_{\cW}')$ is transitive on $\cG$. A $2$--code of $\cS_{3, q}$ is a coclique of $\Gamma_{\cW}'$. We want to apply the Cvetkovi{\'c} bound (Lemma \ref{cvetkovic}) to the graph $\Gamma_{\cW}'$. In order to do that we need to compute the spectrum of $\Gamma_{\cW}'$. Consider the equitable partition arising from the action of the group $\cG_{\Pi_2}$ on $\cG$. Then, according to Lemma \ref{orbits}, the set $\cG$ is partitioned into $\{\{\Pi_2\}, \cG_1, \cG_2, \cG_3, \cG_4\}$ if $q$ is even or into $\{\{\Pi_2\}, \cG_1, \cG_2, \cG_3, \cG_4, \cG_5, \cG_6\}$ if $q$ is odd. Let $B = (b_{i j})$ denote the quotient matrix of this equitable partition. In other words, $b_{i j}$ is the number of planes of $\cG_j$ intersecting a given plane of $\cG_i$ in a line.

\begin{lemma}\label{quotient}
If $q$ is even, then 
$$
B = \begin{pmatrix}
0 & q^3-1 & 0 & 0 & 0 \\
1 & q-2 & 0 & q^3-q & 0 \\
0 & 0 & 0 & q^2-1 & q^2(q-1) \\
0 & q & 1 & q^2-q-2 & q^2(q-1) \\
0 & 0 & 1 & q^2-1 & q^3-q^2-1 \\   
\end{pmatrix}.
$$
If $q$ is odd, then 
$$
B = \begin{pmatrix}
0 & \frac{q^3-1}{2} & \frac{q^3-1}{2} & 0 & 0 & 0 & 0 \\
1 & \frac{q-3}{2} & \frac{q-1}{2} & \frac{q^3-q}{2} & \frac{q^3-q}{2} & 0 & 0 \\
1 & \frac{q-1}{2} & \frac{q-3}{2} & \frac{q^3-q}{2} & \frac{q^3-q}{2} & 0 & 0 \\
0 & \frac{q+1}{2} & \frac{q+1}{2} & \frac{(q-3)(q+1)}{2} & \frac{q^2-1}{2} & \frac{q^2(q-1)}{2} & \frac{q^2(q-1)}{2} \\
0 & \frac{q-1}{2} & \frac{q-1}{2} & \frac{(q-1)^2}{2} & \frac{q^2-1}{2} & \frac{q^2(q-1)}{2} & \frac{q^2(q-1)}{2} \\
0 & 0 & 0 & \frac{q(q-1)}{2} & \frac{q(q+1)}{2} & \frac{q^2(q-1)}{2}-1 & \frac{q^2(q-1)}{2} \\
0 & 0 & 0 & \frac{q(q-1)}{2} & \frac{q(q+1)}{2} & \frac{q^2(q-1)}{2} & \frac{q^2(q-1)}{2}-1 \\
\end{pmatrix}.
$$
\end{lemma}
\begin{proof}
Let $q$ be even. Every plane of $\cG_1$ meets $\Pi_2$ in a line, no plane of $\cup_{i= 2}^4 \cG_{i}$ meets $\Pi_2$ in a line and no plane of $\cG_1$ meets a plane of $\cG_4$ in a line. Hence $b_{1 2} = q^3-1$, $b_{2 1} = 1$, $b_{2 5} = b_{5 2} = b_{1 j} = b_{i 1} = 0$, $j \ne 2$, $i \ne 2$. 

Let $\sigma \in \cG_{1}$ and let $\gamma \in \cup_{i=1}^3 \cG_i$ such that $\gamma \cap \sigma$ is a line. If $\gamma \in \cG_1$, then $\sigma \cap \Pi_2 = \gamma \cap \Pi_2$. Hence $b_{2 2} = q-2$. From Lemma \ref{orbits}, $\gamma \notin \cG_2$ and hence $b_{2 3} = 0$. If $\gamma \in \cG_3$, let $P$ be the point $\gamma \cap \Pi_2$. Then $P \in \sigma \cap \Pi_2$. Let $\ell = \sigma \cap \Sigma_P$ and $s = \gamma \cap \Sigma_P$. Then $|\ell \cap t_P| = 1$ and $|\ell \cap r_P| = 0$. On the other hand the line $s$ is skew to $r_P$ and $t_P$ and meets $\ell$ in a point. Since $s$ can be chosen in $q^2-q$ ways, we have that there are $q^2-q$ planes of $\cG_3$ in $P^\perp$ meeting $\sigma$ in a line. Varying $P$ in $\sigma \cap \Pi_2$, we obtain $b_{2 4} = q^3-q$. 

Let $\sigma \in \cG_{2} \cup \cG_{3}$, $P = \sigma \cap \Pi_2$, and let $\gamma \in \cup_{i=1}^4 \cG_i$ such that $\gamma \cap \sigma$ is a line. Note that necessarily $P \in \gamma$. By Lemma \ref{orbits}, there is no plane of $\cG_2$ meeting a plane of $\cG_1$ in a line and if $\sigma \in \cG_3$, then there are $q$ planes of $\cG_1$ meeting $\sigma$ in a line. Hence $b_{3 2} = 0$ and $b_{4 2} = q$. From Lemma \ref{lemma1}, it follows that $b_{3 3} = 0$, $b_{3 4} = q^2-1$, $b_{4 3} = 1$ and $b_{4 4} = q^2-q-2$. Through a line of $\sigma$ not containing $P$, there pass exactly $q-1$ planes of $\cG_4$. Therefore $b_{4 5} = q^2(q-1)$. 

Let $\sigma \in \cG_4$ and let $\gamma \in \cup_{i=2}^4 \cG_i$ such that $\gamma \cap \sigma$ is a line. If $\gamma \in \cG_2$ or $\cG_3$, then $\sigma \cap \gamma \in \cL_0$ or $\cL_2$. Also, through a line of $\cL_0$ or $\cL_2$ there pass one plane of $\cG_2$ or $\cG_3$ and $q-1$ planes of $\cG_4$, whereas through a line of $\cL_1$ there are $q$ planes of $\cG_4$. Since $\sigma$ contains one line of $\cL_0$, $q+1$ lines of $\cL_1$ and $q^2-1$ lines of $\cL_2$, we have $b_{5 3} = 1$, $b_{5 4} = q^2-1$ and $b_{5 5} = q^3 - q^2 - 1$.

Let $q$ be odd. Every plane of $\cG_1 \cup \cG_2$ meets $\Pi_2$ in a line, no plane of $\cup_{i = 3}^6 \cG_{i}$ meets $\Pi_2$ in a line and no plane of $\cG_1 \cup \cG_2$ meets a plane of $\cG_5 \cup \cG_6$ in a line. Hence $b_{1 2} = b_{1 3} = (q^3-1)/2$, $b_{2 1} = b_{3 1} = 1$, $b_{2 6} = b_{2 7} = b_{3 6} = b_{3 7} = b_{6 2} = b_{7 2} = b_{6 3} = b_{7 3} = b_{1 j} = b_{i 1} = 0$, $j \ne 2, 3$, $i \ne 2, 3$. 

Let $\sigma \in \cG_{1} \cup \cG_{2}$ and let $\gamma \in \cup_{i=1}^4 \cG_i$ such that $\gamma \cap \sigma$ is a line. If $\gamma \in \cG_1 \cup \cG_2$, then $\sigma \cap \Pi_2 = \gamma \cap \Pi_2$. Hence $b_{2 2} = b_{3 3} = (q-3)/2$ and $b_{2 3} = b_{3 2} = (q-1)/2$. If $\gamma \in \cG_3 \cup \cG_4$, let $P$ be the point $\gamma \cap \Pi_2$. Then $P \in \sigma \cap \Pi_2$. The lines $r_P$, $t_P$, $s = \Sigma_P \cap \gamma$ are three pairwise disjoint lines of $\cW_P$. Let $\cR$ be the regulus consisting of the lines of $\cW_P$ intersecting both $r_P$ and $t_P$. The points covered by the lines of $\cR$ form a hyperbolic quadric $\cQ^+(3, q)$. The line $s$ is external or secant to $\cQ^+(3, q)$ according as $\gamma \in \cG_3$ or $\cG_4$. Moreover the line $\ell = \sigma \cap \Sigma_P$ is tangent to $\cQ^+(3, q)$ at the point $\ell \cap t_P$. For a point $L \in \ell \setminus t_P$, the plane $L^\perp \cap \Sigma_P$ meets $\cQ^+(3, q)$ in a non--degenerate conic $C$; in this plane, through the point $L$ there are $(q-1)/2$ lines of $\cW_P$ external to $C$ and skew to both $r_P$, $t_P$ and $(q-1)/2$ lines of $\cW_P$ secant to $C$ and skew to both $r_P$, $t_P$. By varying the point $P$ over the line $\sigma \cap \Pi_2$, we get $b_{2 4} = b_{2 5} = b_{3 4} = b_{3 5} = (q^3 -q)/2$.

Let $\sigma \in \cG_{3} \cup \cG_{4}$, $P = \sigma \cap \Pi_2$, and let $\gamma \in \cup_{i=1}^6 \cG_i$ such that $\gamma \cap \sigma$ is a line. Note that necessarily $P \in \gamma$. As before, let $\cR$ be the regulus consisting of the lines of $\cW_P$ intersecting both $r_P$ and $t_P$ and denote by $\cQ^+(3, q)$ the corresponding hyperbolic quadric. Let $\ell$ be the line $\sigma \cap \Sigma_P$ and $s = \gamma \cap \Sigma_P$. The line $\ell$ is skew to $r_P, t_P$ and it is external or secant to $\cQ^+(3, q)$ according as $\sigma \in \cG_3$ or $\cG_4$. If $\gamma \in \cG_1 \cup \cG_2$, then $s$ is a line of $\cW_P$ meeting both $t_P$ and $\ell$, and it is disjoint from $r_P$. Also $s \cap \ell$ belongs to $\cP_1$ or $\cP_2$, according as $\gamma \in \cG_1$ or $\cG_2$, respectively. From the proof of Lemma \ref{lemma2}, we have that $|\ell \cap \cP_1| = |\ell \cap \cP_2| = (q-1)/2$ if $\ell$ is secant and $|\ell \cap \cP_1| = |\ell \cap \cP_2| = (q+1)/2$ if $\ell$ is external. Hence $b_{4 2} = b_{4 3} = (q+1)/2$ and $b_{5 2} = b_{5 3} = (q-1)/2$. If $\gamma \in \cG_3 \cup \cG_4$, then $s$ is a line of $\cW_P$ intersecting $\ell$ and disjoint from $r_P$ and $t_P$. Also $|s \cap \cQ^+(3, q)|$ equals $0$ or $2$ (i.e., $\ell$ belongs to $\cL_0$ or $\cL_2$, respectively), according as $\gamma \in \cG_3$ or $\cG_4$. If $\ell$ is external, through a point of $\ell$, there are $(q-1)/2$ lines of $\cW_P$ secant to $\cQ^+(3, q)$ and skew to $r_P$ and $t_P$ and $(q-3)/2$ lines of $\cW_P$ external to $\cQ^+(3, q)$ distinct from $\ell$ and skew to $r_P$ and $t_P$. Hence $b_{4 4} = (q-3)(q+1)/2$ and $b_{4 5} = (q^2-1)/2$. If $\ell$ is secant, through a point of $\ell$ not on $\cQ^+(3, q)$, there are $(q-3)/2$ lines of $\cW_P$ secant to $\cQ^+(3, q)$ distinct from $\ell$ and skew to $r_P$ and $t_P$ and $(q-1)/2$ lines of $\cW_P$ external to $\cQ^+(3, q)$ and skew to $r_P$ and $t_P$; through a point of $\ell \cap \cQ^+(3, q)$, there are $q- 1$ lines of $\cW_P$ secant to $\cQ^+(3, q)$ distinct from $\ell$ and skew to $r_P$ and $t_P$. Hence $b_{5 4} = (q-1)^2/2$ and $b_{5 5} = (q^2-1)/2$. A line of $\sigma$ not containing $P$ belongs to $\cL_0$ or $\cL_2$ and, by Lemma \ref{lemma2}, through such a line there pass exactly $(q-1)/2$ planes of $\cG_5$ and $(q-1)/2$ planes of $\cG_6$. Therefore $b_{4 6} = b_{4 7} = b_{5 6} = b_{5 7} = q^2(q-1)/2$. 

Let $\sigma \in \cG_5 \cup \cG_6$ and let $\gamma \in \cup_{i=3}^6 \cG_i$ such that $\gamma \cap \sigma$ is a line. If $\gamma \in \cG_3$ or $\cG_4$, then $\sigma \cap \gamma \in \cL_0$ or $\cL_2$, respectively. Also, through a line of $\cL_0$ or $\cL_2$ there pass one plane of $\cG_3$ or $\cG_4$, $(q-1)/2$ planes of $\cG_5$ and $(q-1)/2$ planes of $\cG_6$, see Lemma \ref{lemma2}. Moreover, by Lemma \ref{lemma2}, through a line of $\cL_1$ there are $q$ planes disjoint from $\Pi_1$ and $\Pi_2$ and they belong either to $\cG_5$ or to $\cG_6$ according as $\sigma\in\cG_5$ or $\sigma\in\cG_6$, respectively. Since, by Corollary \ref{corollary}, $\sigma$ contains $q(q-1)/2$ lines of $\cL_0$, $q+1$ lines of $\cL_1$ and $q(q+1)/2$ lines of $\cL_2$, it follows that $b_{6 4} = b_{7 4} = q(q-1)/2$, $b_{6 5} = b_{7 5} = q(q+1)/2$, $b_{6 6} = b_{77} = q^2(q-1)/2-1$ and $b_{6 7} = b_{7 6} = q^2(q - 1)/2$.     
\end{proof}

\begin{theorem}
The spectrum of the graph $\Gamma_{\cW}'$ is
$$
(q^3-1)^1, (q^2-1)^{\frac{q(q+1)(q^3-1)}{2}}, (-1)^{(q^3-q^2+1)(q^3-1)}, (-q^2-1)^{\frac{q(q-1)(q^3-1)}{2}}.
$$
\end{theorem}
\begin{proof}
The matrix $B$ described in Lemma \ref{quotient} has four distinct eigenvalues, three of them are simple: $q^3-1$, $q^2-1$ and $-q^2-1$, whereas the multiplicity of the eigenvalue $-1$ is two or four, according as $q$ is even or odd. From Lemma \ref{equitable}, the graph $\Gamma_{\cW}'$ has four distinct eigenvalues: $q^3-1$, $q^2-1$, $-1$, $-q^2-1$ with multiplicities $m_0, m_1, m_2, m_3$, respectively. Note that $m_0 = 1$, since $\Gamma_{\cW}'$ is connected. Moreover, the following equations have to be satisfied (see for instance \cite[p. 142]{van}):
\begin{align*}
1 + m_1 + m_2 + m_3 & = q^6, \\
q^3-1 + (q^2-1) m_1 - m_2 - (q^2+1) m_3 & = 0, \\
(q^3-1)^2 + (q^2-1)^2 m_1 + m_2 + (q^2+1)^2 m_3 & = q^6(q^3-1). \\
\end{align*}
It follows that $m_1 = q(q+1)(q^3-1)/2$, $m_2 = (q^3-1)(q^3-q^2+1)$, $m_3 = q(q-1)(q^3-1)/2$.
\end{proof} 

By applying the Cvetkovi{\'c} bound (Lemma \ref{cvetkovic}), we get $\alpha(\Gamma_{\cW}') \le \frac{q(q^2-1)(q^2+q+1)}{2} + 1$. 

\begin{cor}\label{upper_sym3}
Let $\cC$ be a $2$--code of $\cS_{3, q}$, then $|\cC| \le \frac{q(q^2-1)(q^2+q+1)}{2} + 1$.
\end{cor}

\begin{prob}
We obtained a better upper bound for a $2$--code of $\cS_{n, q}$ in the case $n = 3$, by applying the Cvetkovi{\'c} bound to the graph $\Gamma_{\cW}'$, the last subconstituent of $\Gamma_{\cW}$. Determine whether or not this holds true for $n > 3$.
\end{prob}

\subsubsection{The lower bound}


Here we provide the first infinite family of $2$--codes of $\cS_{3, q}$ whose size is larger than the largest possible additive $2$--code. 

\begin{cons}\label{code1}
Let $\cF_{P}$ be a line--spread of $\cW_P$ containing $r_P$ and $t_P$ and let $\cX_{P}$ be the set of $q^2-1$ generators of $\cW(5, q)$ passing through $P$ and meeting $\Sigma_P$ in a line of $\cF_{P} \setminus \{r_P, t_P\}$. Define the set $\cX$ as follows
$$
\bigcup_{P \in \Pi_2} \cX_{P} \cup \{\Pi_2\}.
$$ 
\end{cons}

\begin{theorem}
The set $\cX$ consists of $(q+1)(q^3-1)+1$ planes of $\cW(5, q)$ disjoint from $\Pi_1$ and pairwise intersecting in at most one point.
\end{theorem}
\begin{proof}
By construction every member of $\cX$ distinct from $\Pi_2$ meets $\Pi_2$ in exactly one point. Let $\sigma_1, \sigma_2 \in \cX \setminus \{\Pi_2\}$. If $\sigma_1 \cap \Pi_2 = \sigma_2 \cap \Pi_2$, then $|\sigma_1 \cap \sigma_2| = 1$ and there is nothing to prove. Hence let $P_i = \sigma_i \cap \Pi_2$, $i = 1,2$, with $P_1 \ne P_2$. Assume by contradiction that $\sigma_1 \cap \sigma_2$ is a line, say $\ell$. If $\ell \cap \Pi_2$ is a point, say $R$, then $R \in \sigma_1 \cap \sigma_2$. If $R \ne P_1$ then the line $R P_1$ would be contained in $\sigma_1 \cap \Pi_2$, contradicting the  fact that $|\sigma_1 \cap \Pi_2| = 1$. Similarly if $R \ne P_2$, then $R P_2 \subseteq \sigma_2 \cap \Pi_2$, a contradiction. Therefore $R = P_1 = P_2$, contradicting the fact that $P_1 \ne P_2$. Hence $|\ell \cap \Pi_2| = 0$ and both $\sigma_1, \sigma_2$ are contained in $\ell^\perp$. However in this case the line $P_1 P_2$ is a line of $\cW(5, q)$ contained in $\ell^\perp$ and disjoint from $\ell$; a contradiction. 
\end{proof}

\begin{cor}
There exists a $2$--code of $\cS_{3, q}$ of size $(q^2-1)(q^2+q+1)+1$.
\end{cor}

From Corollary \ref{upper_sym3} a $2$--code of $\cS_{3, 2}$ has at most $22$ elements and hence the $2$--code of $\cS_{3, 2}$ obtained from Construction \ref{code1} is maximal; an alternative proof of its maximality will be exhibited (Corollary \ref{maximal}). Moreover, from \cite{S}, this code is the unique largest $2$--code of $\cS_{3, 2}$ of size $22$. Our next step is to show that, if $q > 2$, the $2$--code of $\cS_{3, q}$ provided in Construction \ref{code1} can be further enlarged. In order to do that some preliminary results are required.

\begin{lemma}\label{lemma_1}
Let $\sigma \notin \cX$ be a plane of $\cW(5, q)$ disjoint from $\Pi_1$ and meeting $\Pi_2$ in at least one point. Then there exists a plane of $\cX$ meeting $\sigma$ in a line.
\end{lemma}
\begin{proof}
Let $\sigma \notin \cX$ be a plane disjoint from $\Pi_1$ and meeting $\Pi_2$ in one point, say $P$. Then $\sigma$ meets $\Sigma_P$ in a line, say $s$, and there are $q+1$ lines of $\cF_P \setminus \{r_P, t_P\}$ meeting $s$ in one point. It follows that there are $q+1$ planes of $\cX_P$ meeting $\sigma$ in a line.
\end{proof}

\begin{lemma}
Through a point $R$ of $\cW(5, q) \setminus (\Pi_1 \cup \Pi_2)$ there pass either $q$ or $q+1$ planes of $\cX$, according as the line through $R$ intersecting $\Pi_1$ and $\Pi_2$ is a line of $\cW(5, q)$ or it is not. 
\end{lemma}
\begin{proof}
Let $R$ be a point of $\cW(5, q) \setminus (\Pi_1 \cup \Pi_2)$, let $\ell_R$ be the unique line through $R$ intersecting both $\Pi_1$, $\Pi_2$ and let $R_i = \ell_R \cap \Pi_i$, $i = 1,2$. Let $s$ denote the line $R^\perp \cap \Pi_2$. If a plane of $\cX_P$ contains the point $R$, then $P \in s$. On the other hand for a fixed $P \in s$ there is at most one plane of $\cX_P$ containing $R$. Hence there are at most $q+1$ planes of $\cX$ through $R$. If $P \in s$ and $P \ne R_2$, then both $R_2$ and $R$ are in $P^\perp$. Hence $\ell_R \subseteq P^\perp$, the line $r_P = P^\perp \cap \Pi_1$ contains $R_1$ and $\langle P, r_P \rangle \cap \ell_R = \{R_1\}$. Therefore $R \notin \langle P, r_P \rangle$ and there exists a plane of $\cX_P$ containing $R$. On the other hand, if $P = R_2$, then $R \in \langle P, r_P \rangle$. Finally note that $R_2 \in s$ if and only if $\ell_R$ is a line of $\cW(5, q)$. 
\end{proof}

Let $\cL$ be the set of lines of $\cW(5, q)$ disjoint from $\Pi_1 \cup \Pi_2$ contained in a plane of $\cX$. Then $|\cL| = q^2(q+1)(q^3-1)$. 

\begin{lemma}\label{rette}
If $q$ is even, then $\cL \subseteq \cL_2$. If $q$ is odd, then $|\cL \cap \cL_0| = |\cL \cap \cL_2|$. 
\end{lemma}
\begin{proof}
Let $\ell$ be a line of $\cL$. Thus there is a point $P \in \Pi_2$ and a plane $\sigma$ of $\cX_P$ containing $\ell$. The three--space $T_\ell$ is contained in $P^\perp$ and does not contain $P$, otherwise $|\sigma \cap \Pi_1| \ne 0$. Hence $T_\ell \cap \cW(5, q)$ is a non--degenerate symplectic polar space $\cW(3, q)$ and $|\cL \cap \cL_1| = 0$. Note that $\cD = \{T_\ell \cap \gamma \; \colon \gamma \in \cX_P\} \cup \{r_\ell, t_\ell\}$ is a line--spread of $\cW(3, q)$. In the point--line dual of $\cW(3, q)$, the line--spread $\cD$ is an ovoid $\cO$ of the parabolic quadric $\cQ(4, q)$ and the regulus determined by $r_\ell, t_\ell, \ell$, would correspond to three points $P_1, P_2, P_3$ of a conic $C$ of $\cQ(4, q)$ such that $P_1, P_2, P_3 \in \cO$.

Assume that $q$ is even. Then the parabolic quadric $\cQ(4, q)$ has a nucleus $N$. If $\ell$ were in $\cL_0$, then $N \in \langle C \rangle$. Consider a three--space $Z$ of the ambient projective space of $\cQ(4, q)$ such that $N \notin Z$. By projecting points and lines of $\cQ(4, q)$ from $N$ to $Z$, we obtain the points and lines of a non--degenerate symplectic polar space $\cW$ of $Z$. In particular $C' = \{NP \cap Z \; \colon P \in C\}$ is a line of $Z$ and $\cO' = \{NP \cap Z \; \colon P \in \cO\}$ is an ovoid of $\cW$. Then we would have $|C' \cap \cO'| \ge 3$, a contradiction, see \cite{Thas}.      

Assume that $q$ is odd. The line $\ell$ belongs to $\cL_2$ if and only if the line polar to the plane $\langle P_1, P_2, P_3 \rangle$ with respect to the polarity of $\cQ(4, q)$ is secant to $\cQ(4, q)$. Let $A$ be the conic obtained by intersecting $\cQ(4, q)$ with the plane polar to the line $\langle P_1, P_2 \rangle$ with respect to the orthogonal polarity of $\PG(4, q)$ associated with $\cQ(4, q)$. Let us count the triple $(R, S, P_3)$, where $R, S \in A$, $R \ne S$, $P_3 \in \cO \setminus \{P_1, P_2\}$ and both $R P_3$, $S P_3$ are lines of $\cQ(4, q)$. The point $R$ can be chosen in $q+1$ ways and for a fixed $R$, the point $P_3$ can be chosen in $q-1$ ways. Finally once $R$ and $P_3$ are fixed, the point $S$ is uniquely determined. Hence there are $q^2-1$ such triples. It turns out that there are $(q^2-1)/2$ points $P_3 \in \cO \setminus \{P_1, P_2\}$ such that the line polar to the plane $\langle P_1, P_2, P_3 \rangle$ with respect to the polarity of $\cQ(4, q)$ is secant to $\cQ(4, q)$ and $(q^2-1)/2$ points $P_3 \in \cO \setminus \{P_1, P_2\}$ such that the line polar to the plane $\langle P_1, P_2, P_3 \rangle$ with respect to the polarity of $\cQ(4, q)$ is external to $\cQ(4, q)$.   
\end{proof}

\begin{cor}\label{maximal}
The $2$--code of $\cS_{3, 2}$ obtained from Construction \ref{code1} is maximal.
\end{cor}
\begin{proof}
From Lemma \ref{lemma_1}, if there exists a plane $\sigma$ disjoint from $\Pi_1$ such that it meets every plane of $\cX$ in at most one point, then $\sigma$ must be disjoint from $\Pi_2$. From Lemma \ref{segre_symp}, $\sigma$ contains exactly one line of $\cL_0$, $3$ lines of $\cL_1$ and $3$ lines of $\cL_2$. Since $|\cL_2| = |\cL|$, the result follows.  
\end{proof}

Let $\Pi_3$ be a generator of $\cW(5, q)$ disjoint from both $\Pi_1$ and $\Pi_2$. Let us denote by $\Pi_i$, $1 \le i \le q+1$, the $q+1$ planes of the unique symplectic Segre variety of $\cW(5, q)$ containing $\Pi_1, \Pi_2, \Pi_3$. In what follows we want to prove that it is possible to construct $\cX$ in such a way that the $q-1$ planes $\Pi_i$, $3 \le i \le q+1$ can be added to it.

\paragraph{The even characteristic case}

Assume that $q > 2$ is even. Since the planes $\Pi_1, \dots, \Pi_{q+1}$ are pairwise disjoint generators of $\cW(5, q)$, from Lemma \ref{segre_symp}, there is a non--degenerate pseudo--polarity $\rho_i$ of $\Pi_i$. The set of absolute points of $\rho_i$ are those of a line $v_i$ of $\Pi_i$. Let $V_i = v_i^{\rho_i}$. Note that the unique line of $\cL_0$ contained in $\Pi_i$ is $v_i$, $3 \le i \le q+1$, while the $q+1$ lines of $\cL_1$ contained in $\Pi_i$ are those through $V_i$, $3 \le i \le q+1$. 

Let $Q$ be a point of $\Pi_2$ not on $v_2$ and distinct from $V_2$. Let $\Sigma_{Q}$ be a $3$--space contained in $Q^\perp$ and not containing $Q$. In particular we choose $\Sigma_{Q}$ spanned by the lines $Q^\perp \cap \Pi_1$ and $Q^{\rho_2}$. Note that $\Sigma_{Q} \cap \Pi_i = Q^\perp \cap \Pi_i$. Indeed, if $s$ is the unique line (not of $\cW(5, q)$) containing $Q$ and meeting each of the planes $\Pi_i$, $1 \le i \le q+1$, in one point, then $(s \cap \Pi_i)^{\rho_i} = \Sigma_Q \cap \Pi_i$. When restricted on $\Sigma_Q$, the polarity $\perp$ defines a non--degenerate symplectic polar space of $\Sigma_Q$, say $\cW_{Q}$. As before, let $r_Q = \Sigma_Q \cap \Pi_1$ and $Q^{\rho_2} = t_Q = \Sigma_Q \cap \Pi_2$. Let $\cR_Q$ be the set of $q+1$ lines of $\cW_Q$ defined as follows
$$
\{\Sigma_Q \cap \Pi_i \; \colon 1 \le i \le q+1\}.
$$
Then $\cR_Q$ is a regulus of $\cW_Q$ containing both $r_Q$ and $t_Q$; the opposite regulus of $\cR_Q$ contains exactly one line of $\cW_Q$ which is the line consisting of the points $v_i\cap(s\cap\Pi_i)^{\rho_i}$. 

\begin{lemma}
There exists a Desarguesian line--spread of $\cW_P$ having in common with $\cR_Q$ exactly the lines $r_P$ and $t_P$.
\end{lemma}
\begin{proof}
Let $\cQ(4, q)$ be the point line dual of $\cW_Q$ and let $N$ be the nucleus of $\cQ(4, q)$. The regulus $\cR_Q$ corresponds to a conic $C$ of $\cQ(4, q)$ such that $N \notin \langle C \rangle$ and  the lines $r_P$ and $t_P$ correspond to two points of $C$, say $R$ and $T$. The result follows, since there are $q^2/2-q$ elliptic quadrics of $\cQ(4, q)$ meeting $C$ exactly in the points $R, T$.  
\end{proof}

For any point $Q$ of $\Pi_2$ different from $V_2$ and not on $v_2$, let $\cF_{Q}$ be a Desarguesian line--spread of $\cW_Q$ having in common with $\cR_Q$ exactly the lines $r_Q$ and $t_Q$ and let $\cY_{Q}$ be the set of $q^2-1$ generators of $\cW(5, q)$ passing through $Q$ and meeting $\Sigma_Q$ in a line of $\cF_{Q} \setminus \{r_Q, t_Q\}$. For any point $P \in v_2 \cup \{V_2\}$, let $\cX_{P}$ be a set of $q^2-1$ generators of $\cW(5, q)$ passing through $P$ defined as in Construction \ref{code1}. 

Define the set ${\bar \cX}$ as follows
$$
\left( \bigcup_{P \in v_2 \cup \{V_2\}} \cX_{P} \right) \cup \left( \bigcup_{Q \in \Pi_2 \setminus (v_2 \cup \{V_2\})} \cY_{Q} \right) \cup \left( \bigcup_{i = 2}^{q+1} \Pi_i \right).
$$ 
\begin{theorem}
The set $\bar{\cX}$ consists of $q^4 + q^3 +1$ planes of $\cW(5, q)$ disjoint from $\Pi_1$ and pairwise intersecting in at most one point.
\end{theorem}
\begin{proof}
It is enough to show that a plane $\sigma$ of 
$$
\left( \bigcup_{P \in v_2 \cup \{V_2\}} \cX_{P} \right) \cup \left( \bigcup_{Q \in \Pi_2 \setminus (v_2 \cup \{V_2\})} \cY_{Q} \right) 
$$
meets $\Pi_i$, $3 \le i \le q+1$, in at most one point. If $\sigma$ intersects $\Pi_i$ in a line, say $\ell$, from Lemma~\ref{rette}, $\ell \in \cL \subseteq \cL_2$. Hence $\ell \ne v_i$ and $V_i \notin \ell$. Let $s$ be the unique line through the point $\ell^{\rho_i}$ meeting both $\Pi_1$ and $\Pi_2$ in one point. Let $Q = s \cap \Pi_2$. Then $Q$ coincides with $\ell^\perp \cap \Pi_2$ and $Q \notin v_2 \cup \{V_2\}$. Moreover $\sigma = \langle Q, \ell \rangle \in \cY_Q$. But in this case the Desarguesian line--spread $\cF_Q$ would have the three lines $\ell, r_Q, t_Q$ in common with $\cR_Q$, contradicting the fact that $|\cF_Q \cap \cR_Q| = 2$. 
\end{proof}

\paragraph{The odd characteristic case}

Assume that $q$ is odd. Since the planes $\Pi_1, \dots, \Pi_{q+1}$ are pairwise disjoint generators of $\cW(5, q)$, from Lemma \ref{segre_symp}, there is a non--degenerate orthogonal polarity $\rho_i$ of $\Pi_i$. The set of absolute points of $\rho_i$ are those of a conic $\alpha_i$ of $\Pi_i$. Note that a line $\ell$ of $\Pi_i$ belongs to $\cL_j$, according as $|\ell \cap \alpha_i| = j$, $0 \le j \le 2$, $3 \le i \le q+1$. 

Let $Q$ be a point of $\Pi_2$ not on $\alpha_2$. Let $\Sigma_{Q}$ be a $3$--space contained in $Q^\perp$ and not containing $Q$. In particular we choose $\Sigma_{Q}$ spanned by the lines $Q^\perp \cap \Pi_1$ and $Q^{\rho_2}$. Note that $\Sigma_{Q} \cap \Pi_i = Q^\perp \cap \Pi_i$. Indeed, if $s$ is the unique line containing $Q$ and meeting each of the planes $\Pi_i$, $1 \le i \le q+1$, in one point, then $(s \cap \Pi_i)^{\rho_i} = \Sigma_Q \cap \Pi_i$. When restricted on $\Sigma_Q$, the polarity $\perp$ defines a non--degenerate symplectic polar space of $\Sigma_Q$, say $\cW_{Q}$. As before, let $r_Q = \Sigma_Q \cap \Pi_1$ and $t_Q = \Sigma_Q \cap \Pi_2$. Let $\cR_Q$ be the set of $q+1$ lines of $\cW_Q$ defined as follows
$$
\{\Sigma_Q \cap \Pi_i \; \colon 1 \le i \le q+1\}.
$$
Then $\cR_Q$ is a regulus of $\cW_Q$ containing both $r_Q$ and $t_Q$ and the opposite regulus of $\cR_Q$ contains exactly $0$ or $2$ lines of $\cW_Q$. The proof of the next result is left to the reader.

\begin{lemma}
There exists a Desarguesian line--spread of $\cW_P$ having in common with $\cR_Q$ exactly the lines $r_P$ and $t_P$.
\end{lemma}

For any point $Q$ of $\Pi_2$  not on $\alpha_2$, let $\cF_{Q}$ be a Desarguesian line--spread of $\cW_Q$ having in common with $\cR_Q$ exactly the lines $r_Q$ and $t_Q$ and let $\cY_{Q}$ be the set of $q^2-1$ generators of $\cW(5, q)$ passing through $Q$ and meeting $\Sigma_Q$ in a line of $\cF_{Q} \setminus \{r_Q, t_Q\}$. For any point $P \in \alpha_2$, let $\cX_{P}$ be a set of $q^2-1$ generators of $\cW(5, q)$ passing through $P$ defined as in Construction \ref{code1}. 

Define the set ${\bar \cX}$ as follows
$$
\left( \bigcup_{P \in \alpha_2} \cX_{P} \right) \cup \left( \bigcup_{Q \in \Pi_2 \setminus \alpha_2} \cY_{Q} \right) \cup \left( \bigcup_{i = 2}^{q+1} \Pi_i \right).
$$ 

A proof similar to that given in the even characteristic case yields the following result.

\begin{theorem}
The set $\bar{\cX}$ consists of $q^4 + q^3 +1$ planes of $\cW(5, q)$ disjoint from $\Pi_1$ and pairwise intersecting in at most one point.
\end{theorem}

\subsection{$2$--codes of $\cH_{3, q^2}$}

Let $\bar{\perp}$ be the Hermitian polarity of $\PG(5, q^2)$ defining $\cH(5, q^2)$. Recall that $\bar{\cG}$ is the set of $q^9$ planes of $\cH(5, q^2)$ that are disjoint from $\Lambda_1$, the group $\bar{G}$ is the stabilizer of $\Lambda_1$ in $\PGU(6, q^2)$, $\Lambda_2 = L(0_3)$, and $\bar{G}_{\Lambda_2}$ is the stabilizer of $\Lambda_2$ in $\bar{G}$. For a point $P$ in $\Lambda_2$, let $\bar{\Sigma}_P$ denote a $3$--space contained in $P^{\bar{\perp}}$ and not containing $P$. When restricted to $\bar{\Sigma}_P$, the polarity $\bar{\perp}$ defines a non--degenerate Hermitian polar space of $\bar{\Sigma}_P$, say $\cH_{P}$. Moreover $\bar{r}_P = \bar{\Sigma}_P \cap \Lambda_1$ and $\bar{t}_P = \bar{\Sigma}_P \cap \Lambda_2$ are lines of $\cH_P$. 

\begin{cons}\label{code2}
Let $\bar{\cF}_{P}$ be a partial spread of $\cH_P$ containing $\bar{r}_P$ and $\bar{t}_P$ and let $\cY_{P}$ be the set of $|\bar{\cF}_P| - 2$ generators of $\cH(5, q^2)$ passing through $P$ and meeting $\bar{\Sigma}_P$ in a line of $\bar{\cF}_{P} \setminus \{\bar{r}_P, \bar{t}_P\}$. Define the set $\cY$ as follows
$$
\bigcup_{P \in \Lambda_2} \cY_{P} \cup \{\Lambda_2\}.
$$ 
\end{cons}

A proof similar to that given in the symmetric case gives:

\begin{theorem}
The set $\cY$ consists of $(q^4+q^2+1)(|\cF| - 2)$ planes of $\cH(5, q^2)$ disjoint from $\Lambda_1$ and pairwise intersecting in at most one point.
\end{theorem}

By selecting $\cF$ as a partial spread of $\cH(3, q^2)$ of size $(3q^2-q+2)/2$, see \cite[p. 32]{ACE}, the following arises. 

\begin{cor}
There exists a $2$--code of $\cH_{3, q^2}$ of size $q^6 + \frac{q(q-1)(q^4+q^2+1)}{2}$.
\end{cor}

\smallskip
{\footnotesize
\noindent\textit{Acknowledgments.}
This work was supported by the Italian National Group for Algebraic and Geometric Structures and their Applications (GNSAGA-- INdAM).}

\end{document}